\documentclass[11pt]{amsart}
\usepackage{amscd}
\usepackage{amsmath,amssymb,latexsym}
\usepackage{pb-diagram}
\usepackage{dynkin-diagrams}
\usepackage{enumerate}
\usepackage{amsfonts}
\usepackage{amsthm}
\usepackage{mathtools}
\usepackage{tikz}\usetikzlibrary{chains}
\usepackage{hyperref}

\newcommand{\be}{\begin{enumerate}}
\newcommand{\ee}{\end{enumerate}}
\newcommand{\beq}{\begin{equation}}
\newcommand{\eeq}{\end{equation}}
\newcommand{\integral}[1]{\int\limits_{ {#1}(F)\backslash {#1}(\A)}}

\newtheorem{Thm}{Theorem}[section]
\newtheorem{Cor}{Corollary}[section]
\newtheorem{Lem}{Lemma}[section]
\newtheorem{Prop}{Proposition}[section]
\newtheorem{Claim}{Claim}[section]
\newtheorem{Remark}[equation]{Remark}
\newtheorem{Def}[equation]{Definition}
\newcommand{\ra}{\rangle}
\newcommand{\la}{\langle}

\newcommand{\C}{\mathbb C}

\newcommand{\A}{\mathbb A}
\newcommand{\bF}{\mathbb F}

\newcommand{\mB}{\mathcal B}

\newcommand{\mS}{\mathcal S}
\newcommand{\mE}{\mathcal E}

\newcommand{\Z}{\mathbb Z}

\newcommand{\mM}{\mathcal M}
\newcommand{\mO}{\mathcal O}
\newcommand{\mP}{\mathcal P}

\newcommand{\G}{\mathbb G}
\newcommand{\mg}{\mathfrak g}
\newcommand{\mA}{\mathcal A}

\newcommand\Ker{\operatorname{Ker}}

\newcommand\ind{\operatorname{ind}}
\newcommand\Ind{\operatorname{Ind}}

\newcommand\Hom{\operatorname{Hom}}
\newcommand\cha{\operatorname{char}}
\newcommand\Tr{\operatorname{Tr}}
\newcommand\Res{\operatorname{Res}}

 \newcommand{\sms}{\smallskip}

\usepackage{xcolor}
\usepackage{amscd}

\theoremstyle{plain}

\tikzset
    {node distance=2em, ch/.style=
      {circle,draw,on chain,inner sep=2pt},chj/.style=
      {ch,join},every path/.style={shorten >=4pt,shorten <=4pt},
      line width=1pt,baseline=-1ex}

\title[Automorphic functionals
  for the minimal representations]{Automorphic functionals
  for the minimal representations of groups of type
  $D_n$ and $E_n$}
\author{Nadya Gurevich and David Kazhdan}

\begin{document}
\maketitle

\hskip 3cm{\sl{In the memory of an outstanding mathematician}}

\hskip 3cm{\sl{and a remarkable person Anatoly  Vershik}}

\begin{abstract}
Let $G$ be split simply-connected group of type $D$ or $E$.
The minimal automorphic representation $\Pi$ of $G(\A)$
admits a realization on a space of functions $\mS(X(\A))$ for a variety $X$. 
In this paper we write explicitly an automorphic, i.e. $G(F)$-invariant,
functional on $\mS(X(\A)).$
\end{abstract}

\section{Introduction}\label{sec:intro}

\subsection{Groups}

Let $G$ be a split, almost simple, simply-connected Chevalley group
with Lie algebra $\mg$. Such $G$ is uniquely determined by its
root system, which is irreducible. The root system $R$ is either
of a classical type $A_n,B_n,C_n, D_n$ or of one of $5$ exceptional
types: $E_6,E_7,E_8, F_4$ or $G_2$.

Let $F$ be the field of rational functions on
a smooth, projective, absolutely irreducible curve $C$
over a finite field $\bF_q$ of characteristic $p\neq 2$ 
and $\A$ be  the ring of adeles of $F$.
For a closed point $v\in C$  we denote by $F_v$ the local non-archimedian
field, which is the completion of $F$ at $v$.

In this paper  we analyze a realization of 
the minimal automorphic representation $\Pi$ of
the adelic group $G(\A)$ in a space of functions $\mS(X(\A))$
on a variety $X$ and describe explicitly the  $G(F)$-invariant
functional on $\mS(X(\A))$.

One has  $\Pi\simeq\otimes_v \Pi_v,$ where the representations 
$\Pi_v$ of $G(F_v)$ are minimal.

\subsection{Local theory}\label{intro:local:theory}

In this subsection we fix a non-archimedian local field $F$
of characteristic either zero or positive,
with the ring of integers $\mO_F$, maximal ideal $\mP$, generated by a uniformizer $\varpi$ and the residue field of cardinality $q$. 

There exist several equivalent definitions of minimal representation.
We recall the definition
that is convenient for our setting, following the approach
of \cite{GanSavin}.

The group $G$ contains unique, up to conjugation,
parabolic subgroup, whose 
unipotent radical $V$ is a Heisenberg group.
We denote by $Z$ the one-dimensional center of $V$.
This parabolic subgroup  is maximal, unless $G$ is of type $A_n$.

Let $\psi$ be a non-trivial complex character of $Z(F)$. For any
smooth representation $(\pi,\mathcal V_\pi)$ of $G(F)$, the space of $(Z,\psi)$ coinvariants
$$\pi_{Z,\psi}=\mathcal V_\pi/Span\{\pi(z)v-\psi(z)v, z\in Z(F),v\in \mathcal V_\pi\}$$
is a representation of $V(F)$ with central character $\psi$.

If for an irreducible smooth representation
$\pi$ of $G(F)$ one has  $\pi_{Z,\psi}=0$ for all non-trivial $\psi$,
then it is easy to see that $\pi$ is one-dimensional.


\begin{Def}\label{min:def} A smooth irreducible  representation
  $\Pi$ of $G(F)$ is called
  {\it minimal} if there exists a non-trivial character $\psi$ of $Z(F)$
  such that $\Pi_{Z,\psi}$ is a non-zero irreducible representation of $V(F)$. 
\end{Def}

By Stone-von-Neumann theorem there exists unique irreducible representation
$\rho_\psi$ of $V(F)$ with central character $\psi$.
Hence if $\Pi$ is minimal, then  $\Pi_{Z,\psi}\simeq \rho_\psi$
as representations of $V(F)$. 

The group $V/Z$ admits a structure of non-degenerate symplectic space,
that we denote by $W$.
In the seminal paper \cite{Weil}  Weil  showed that if
$\cha(F)\neq 2$, $\rho_\psi$ can be canonically
extended to a representation of
 $Mp(W)\ltimes V(F)$ for a double cover 
 $Mp(W)$   of $Sp(W)(F)$.
As a representation of $Mp(W)$,  $\rho_\psi$
is a sum of two irreducible representations
$\rho_\psi=\rho_\psi^+\oplus \rho_\psi^-$.
Later Weil has expanded the construction of $\rho_\psi$
to the case $\cha(F)=2$,
we exclude this case from our discussion, since we assume that $\cha(F)\neq 2$.

 The complex dual group $\hat G(\C)$ of $G$ contains
unique subregular unipotent orbit $\mO_{sr}$. By Jacobson-Morozov
theorem to such orbit is associated a semisimple conjugacy class,
denoted by $\{t_{sr}\}$. 

The following statement  summarizes the results on the existence and
uniqueness of minimal representations. 

\begin{Thm}
  \be
  \item
  If $G$ is of type $D_n, n\ge 4$ or $E_n$, then there
  exists unique minimal representation $\Pi$ of $G(F)$.
  It is unramified and its Satake parameter is $\{t_{sr}\}$. 
\item If $G$ is of type $A_n$, then the unramified 
  representation $\Pi$ of $G(F)$ with Satake parameter
  $\{t_{sr}\}$ is minimal.
There exist other minimal representations.   
\item If $G$ is of other types, there are no minimal representations
  of $G(F)$.
  \ee
\end{Thm}

\begin{proof}
  Minimal representation for groups of type $D_n,E_n$ have
  been constructed and their Satake parameters computed in \cite{KazhdanSavin}.
  See \cite{Torasso} for uniqueness.  
  A different construction, which includes the group of type $A_n$
  is given in \cite{Savin:invent}. For part $(3)$
 note that $\Pi_{Z,\psi}$ is a representation of $[P,P]$
 which centralizes the character $\psi$. This implies
 existence of a splitting of $[P,P]$ in $Mp(W)$. Such splitting
 does not exist if $G$ is of type $B_n,C_n$ or $F_4$, which implies
 non-existence of minimal representations of $G(F)$. 
 While the splitting of $[P,P]$  exists for the group $G$ of type $G_2$,
 it is proven in \cite{GanSavin} that $G(F)$ does not possess a minimal representation.

\end{proof}

If $G$ is  of type $C_n,F_4, G_2$ or $B_3$, then  minimal representations exist
on certain central covers of $G(F)$. In the paper we shall only consider
the groups of type $D_n$ and $E_n$ and occasionally $A_n$. The unramified
representation $\Pi$ for these groups with Satake parameter $\{t_{sr}\}$
will be referred to as {\sl the minimal representation.}

\subsubsection{Other definitions}
Let us  explain the relation between our notion of minimality
and the minimal  non-zero $G(\bar F)$-orbit
$\mO_{min}$ in $\mg^\ast(\bar F)$. To simplify the presentation we assume 
that $\mO_{min}\cap \mg^\ast(F)$ is a single $G(F)$-orbit. This assumption
holds for the groups of type $A_n,D_n$ for $n>1$ and $E_n$. 

Assume further that $\cha(F)=0$ and denote by
$\exp :\mg(F)\rightarrow G(F)$  the exponential map.
Any distribution on $G(F)$ can be regarded as a distribution on
a neighborhood of $0\in \mg(F)$.
For any coadjoint nilpotent orbit $\mO$ in $\mg^\ast(F)$
let $\mu_{\mO}$ be the distribution given by properly normalized
$G(F)$-invariant measure of $\mO$ and $\hat \mu_{\mO}$ be its
Fourier transform.

By Howe and Harish-Chandra, for any  smooth representation $\pi$ of $G(F)$
there exists a neighborhood $1\in U\subset G(F)$ and
complex numbers $c_{\mO}(\pi)$
such that for the character distribution $\Theta_\pi$ holds:
$$\Theta_\pi=\sum_{\mO} c_\mO(\pi)\hat \mu_\mO$$
as distributions on $\mS_c(\exp^{-1}(U(F)))$.
The sum runs over nilpotent orbits $\mO$ in $\mg^\ast(F)$. 

The relation of the minimal representation to the minimal orbit is given
in the following proposition.

\begin{Prop} (\cite{GanSavin}, Proposition $3.7$)
  An irreducible representation $\Pi$ is minimal if and only if
  $\Theta_\Pi=\hat\mu_{\mO_{min}}+c_0$.
\end{Prop}

Using the homogeneity property of the distributions $\mu_{\mO}$
and the fact that $\dim \pi^{K}$ for any open compact subgroup $K$ 
 equals  the trace of the projection of $\pi$ to $\pi^{K}$
it is shown in \cite{Savin:invent},
section $2$:

\begin{Cor} An irreducible representation $\Pi$ is minimal if and only if
the Gelfand-Kirillov dimension of $\Pi$ is $D=\frac{1}{2}\cdot\dim \mO_{min}$. 
Precisely, for the sequence of principal congruence subgroups $K_n$,   
$\lim_{n\to \infty}\dim(\Pi^{K_n})\cdot q^{-D}=1$.
\end{Cor}

\begin{Remark}
  If $\cha(F)$ is finite, but sufficiently large,
  the exponential map  can be replaced by the  map $e:\mg(F)\rightarrow G(F)$
  defined in \cite{KazhdanVarshavsky}, which satisfies conditions of Theorem
  $3.5.2$ in \cite{DeBacker}. Therefore, the Harish-Chandra asymptotic expansion
  holds if $char (F)\gg 1$.

\end{Remark}


\subsubsection{Models}
We are interested in constructing models of the minimal representation
$\Pi$, i.e. its realization  on some space of locally constant functions
$\mS(X(F))$ on a variety $X$ of dimension $\frac{1}{2}\cdot \dim (\mO_{min}).$
There could be more than one model.
For every model the action of  elements of certain maximal parabolic subgroup
is easy to describe. We shall consider two types of models:
the Heisenberg model, associated  with the parabolic subgroup  whose unipotent radical is a Heisenberg group 
and the Siegel (sometimes also called Schrodinger) model, associated
to a Siegel parabolic subgroup, whose unipotent radical is abelian. 

Kazhdan and Savin used the Heisenberg model to establish the existence
of the minimal representation. Later Savin \cite{Savin:invent}
has established the existence of the Siegel model
for groups of type $A_n, D_n, E_6$ and $E_7$. The group of type $E_8$
does not possess a Siegel parabolic subgroup.
Since the representation $\Pi$ is unramified,
the space $\mS(X(F))$ contains a one-dimensional space of
spherical functions. The formula for the normalized spherical
function in each type of model was obtained in \cite{KazhdanPolishchuk} and
\cite{SavinWoodbury}.

\subsection{Global theory}
In this subsection let $F$ be the field of rational functions on a
smooth projective curve $C$ and $G$ be of type $D_n$ or $E_n$. 
Let $\Pi\simeq \otimes \Pi_v$ be the minimal representation
of the adelic group $G(\A)$.
Having the models $\mS(X(F_v))$ of $\Pi_v$  for all $v$,
we can form their restricted
tensor product with respect to the choice of local spherical
functions. This results in a model $\mS(X(\A))$
for the minimal representation $\Pi$ of  $G(\A)$.
As the minimal representation $\Pi$  is  automorphic,
there exists a  $G(F)$-invariant linear functional on
$\mS(X(\A)).$ It is natural to ask for an explicit formula for it.

\subsubsection{Metaplectic minimal representation}
 In this subsection we 
 present a motivating example, although it is  outside of the scope
 of this paper.

For any symplectic space $W$, consider the Weil representation
$\rho_{\psi_v}$ of $Mp(W)(F_v)$. It is reducible and the representations
$\rho^\pm_{\psi_v}$ are minimal according to Definition \ref{min:def}.

To describe a model of $\rho_{\psi_v}$ we choose
a presentation of $W=W^+\oplus W^-$ as a direct sum of
two isotropic spaces. The $\rho_{\psi_v}$ admits
a realization on the space $\mS(W^+(F_v))$
of locally constant functions of compact support.
The action of the Siegel maximal parabolic subgroup $\bar P(F_v)$
of $Mp(W)(F_v)$ is explicit and an additional element $\sigma$
such that $\sigma(W^{\pm})=W^{\mp}$  acts  by the classical Fourier transform.
If  $\psi_v$ has conductor  $\mO_v$, and the residual
characteristic of $F_v$ is not $2$,
then the characteristic function  $\phi^0_v=1_{W^+(\mO_v)}$ is the
normalized spherical function in this model.

This allows to define a realization of the (highly reducible)
minimal representation
${\rho_\psi=\otimes \rho_{\psi_v}}$ of the adelic group $Mp(W)(\A)$, depending on
an additive character $\psi: F\backslash \A\rightarrow \C$,
on the space  $\mS(W^+(\A))$, that is a restricted tensor product
of local spaces with respect to the choice of spherical functions above.
It is remarkable, that
the representation $\rho_\psi$ is automorphic and 
 appears as a direct summand in the space of square-integrable
automorphic forms $\mathcal A_2(Mp(W))$.
An embedding of $\rho_\psi$ in the space of automorphic
forms is given by $\phi\mapsto F_\phi$, where
$F_\phi(g)=\Theta(g\cdot\phi)$ and $\Theta$ is  the
$G(F)$-invariant functional on $\mS(W^+(\A))$ given by
\beq \label{theta:def:intro}
\Theta(\phi)=\sum_{x\in W^+(F)} \phi(x).
\eeq
For the functional field $F$ the sum for any $\phi$ is finite,
hence the functional is well-defined. 

To check the $G(F)$-invariance it is enough to check the invariance under
$P(F)$, which follows from  explicit formulas of the action
and the invariance under the element $\sigma\in Sp(W)(F)$
which follows from the classical Poisson summation formula. 

\subsubsection{Main Theorem: decomposition of the automorphic functional
  $\theta_{G,L}$}
Let   $\Pi=\otimes_v \Pi_v$ be the minimal representation  of $G(\A),$
where $G$ is of type  $D_n$ or $E_n$.
Consider   a model $\mS(X(\A))$ of $\Pi$ on which a maximal parabolic
subgroup $P=LV$  acts explicitly.
In the paper  the parabolic $P$ will be of one of two types:
either the unipotent radical $V$ of $P$  is Heisenberg group, as in
the beginning of subsection \ref{intro:local:theory},
or $V$ is abelian. 
For the rest of introduction we assume for simplicity that $V$ is abelian. 
The case when $V$ is a Heisenberg group is
 discussed in the body of the paper.

Let $\mA_2(G)$ denote the space
of square-integrable automorphic forms on $G$. The space 
$\Hom_{G(\A)}(\Pi,\mA_2(G))$ is one-dimensional and
we denote by $\bar\theta(\Pi)\subset \mA_2(G)$ 
the image of $\Pi$. One can show that $\bar\theta(\Pi)$ consists of 
residues of degenerate Eisenstein series. We use  the Fourier
 expansion of an automorphic form in $\bar \theta(\Pi)$ along
 $V(F)\backslash V(\A)$, to write a $G(F)$-invariant functional on $\mS(X(\A))$. 
Since $V$ is abelian,  points $x\in X(F)$ parameterize characters $\psi_x$ of
$V(F)\backslash V(\A)$, such that $(V,\psi_x)$-Fourier coefficients
of functions in $\bar\theta(\Pi)$
 are not identically zero.
 In Section \ref{sec:main:thm} we define a  $G(\A)$-equivariant isomorphism
 $\tau_X$ satisfying
 $$\tau_X: \mS(X(\A)) \xrightarrow{\sim} \bar \theta(\Pi), 
 \quad \tau_X(\phi)^{V,\psi_x}(e)=\phi(x), \quad \forall x\in X(F), \phi\in \mS(X(\A)).$$  
The functional $\theta_{G,L}: \phi\mapsto \tau_X(\phi)(e)$  on
$\mS(X(\A))$ is  $G(F)$-invariant.
The goal of the paper is to write this functional explicitly.
 
A natural candidate to consider would be an averaging functional
$$\theta_X(\phi)=\sum_{x\in X(F)}\phi(x).$$
However, this functional  is only $P(F)$-invariant.
To obtain a $G(F)$-invariant functional
$\theta_{G,L}$  we have to add 
certain boundary terms arising from
the constant term of functions in $\bar\theta(\Pi)$ along $V$. 
Our main Theorem \ref{main:thm} claims:
\begin{equation}\label{main:ne:intro}
  \theta_{G,L}=\theta_X+p_X+|\Delta_F|^{1/2} \theta_{G_1,M_1}\circ \mB,
  \end{equation}
where
\begin{itemize}
\item $G_1=[L,L]$ and $\mS(X_1(\A))$
  is a model of the minimal representation $\Pi_1$ of $G_1(\A)$
  on which certain maximal parabolic subgroup $Q_1=M_1V_1$ of $G_1$
acts by simple formula.

\item $\mB\in \Hom_{V(\A)G_1(\A)}(\mS(X(\A)),\mS(X_1(\A))$
  is  a particular non-zero map  in the one-dimensional space.

\item The functional $\theta_{G_1,M_1}$ is a $G_1(F)$-invariant functional
  on $\mS(X_1(\A))$, defined similarly to $\theta_{G,L}$.

\item
  The functional $p_X\in \Hom_{V(\A)G_1(\A)}(\mS(X(\A)),\C)$
  sends any function to its regularized value at the boundary point
  $0$ of $X$,  see equation (\ref{p:X:global}).

\item $\Delta_F$ is the discriminant of the field $F$. 
\end{itemize}

The $P(F)$-invariant functional $\theta_X$ in the expansion above
is called {\sl the main term}
while the functionals $p_X$ and $\theta_{G_1,M_1}\circ \mB$ that are
$V(\A)G_1(F)$ invariant are called {\sl boundary terms.} 

The construction of the boundary map $\mB$ uses the 
structure  of Jacquet module $(\Pi_v)_V$, that is
available only for non-archimedean fields. For this
reason only we choose $F$ to be a function field. We also
add an assumption that the genus $g$  of the curve $C$ is positive.

If  $V$ is abelian the space $\Hom_{P(F)}(\mS(X(\A)),\C)$
is three-dimensional and  the functionals
$\theta_X,p_X$ and $\theta_{G_1,M_1}\circ \mB$
form its basis. The content of Theorem \ref{main:thm} is to specify the coefficients
in the decomposition of $\theta_{G,L}$ with respect to this basis.
If $V$ is Heisenberg the space $\Hom_{P(F)}(\mS(X(\A)),\C)$ is four-dimensional.
In this case our main Theorem is the refinement of the  Theorem $8.2.1$
in \cite{KazhdanPolishchuk}.

\subsubsection{Expansion of $\theta_{G,L}$}
If the pair $(G_1,M_1)$ is admissible, see Definition \ref{admissible:def},
 Theorem \ref{main:thm}
can be applied to the last boundary term,  giving rise to further terms
in the  decomposition of the functional $\theta_{G,L}$.

Let us illustrate such decomposition for  the pair $(D_n,D_{n-1})$,
that is admissible for $n\ge 4$. The variety associated
to the pair $(D_n,D_{n-1})$ is the cone $X_n$ of non-zero isotropic vectors
in a $2n-2$ dimensional quadratic space. For the  affine closure $\bar X_n$ of $X_n$
one has  $\bar X_n(F)=X_n(F)\cup\{0\}$. The regularized value at $0$
defines the functional $p_n$ on $\mS(X_n(\A))$ and there is
$P(\A)$ equivariant map
$\mB_{n}: \mS(X_n(\A))\rightarrow  \mS(X_{n-1}(\A)).$

We write 
$$\mB_n^{(0)}=Id, \quad
\mB_{n}^{(k)}=\mB_{n-(k-1)}\circ \ldots \circ \mB_{n}:
\mS(X_n(\A))\rightarrow  \mS(X_{n-k}(\A))$$


Applying Theorem \ref{main:thm}
several times we obtain an expansion for $\theta_{D_n,D_{n-1}}$
in terms of maps $\mB_{n}^{(k)}$.
\begin{Thm}\label{Dn:automorphic:decomp}  For $n\ge 4$ one has
$\theta_{D_n,D_{n-1}}=$
\beq \label{Dn:decomp}
  \sum^{n-4}_{k=0}  |\Delta_F|^{\frac{k}{2}} \theta_{X_{n-k}}\circ \mB_{n}^{(k)}+
\sum_{k=0}^{n-4} |\Delta_F|^{\frac{k}{2}} p_{n-k}\circ \mB_n^{(k)}+ 
|\Delta_F|^{\frac{n-3}{2}}\theta_{D_3,D_2}\circ \mB_n^{(n-3)}.
\eeq
\end{Thm}
\begin{Remark}
The summation over isotropic cones of decreasing dimension also
appears in the preprint  \cite{Getz}, where the author uses the description
of the minimal representation of $D_n$ as the theta lift
from the trivial representation of $SL_2$. 
\end{Remark}

Theorem \ref{main:thm}
can not be applied to the last term  in the equation  (\ref{Dn:decomp})
since the pair $(D_3,D_2)$ is not admissible.
In Section $7$ we compute the functional $\theta_{D_3,D_2}$ for the group $SL_4$, thus completing
the decomposition of $\theta_{D_n,D_{n-1}}$. 

We seek to write a similar expansion of the automorphic functional
for a group $G$ of type $E$.
Let us illustrate this for $G=E_6$. The cases   $E_7,E_8$ are similar.
Consider the pair $(G,L)=(E_6,D_5)$ and the associated pair
$(G_1,M_1)=(D_5,A_4)$ that is not admissible. Theorem \ref{main:thm} implies
\beq \label{intro:E6:D5}
\theta_{E_6,D_5}=\theta_X+p_X+|\Delta_F|^{1/2}\cdot \theta_{D_5,A_4}\circ \mB,
\eeq
and Theorem \ref{main:thm} can not be applied again to the last term. 
We prove in {Section \ref{sec:transition}}
\beq \label{intro:D5:A4:D4} 
\theta_{D_5,A_4}=|\Delta_F|^{1/2}\cdot\theta_{D_5,D_4}\circ T
\eeq
where $T$ is a transition map between  the models corresponding to
$(D_5,A_4)$ and $(D_5,D_4)$, which is admissible.
Substituting  (\ref{intro:D5:A4:D4}) in  (\ref{intro:E6:D5})
we can continue with the expansion of $\theta_{D_5,D_4}$ discussed above. 
\subsection{Outline of the paper}
 After setting notation in
Section $2$ we review in Section $3$ the construction and properties
of models of minimal representations for groups over local field.
In Section $4$ we use the local models to
construct the model for the adelic minimal representation.
We define several functionals on it, that are ingredients in the main
theorem. Section $5$ is  devoted to the automorphic representation
$\bar\theta(\Pi)$, spanned by residues 
of a degenerate Eisenstein series, which is isomorphic to $\Pi$.
The constant term along $V$ of any function in $\bar\theta(\Pi)$
is a sum of two terms $J_{v_0}$ and $J_{v_1},$
where the elements $v_0, v_1$ in the Weyl group of $G$
are defined for each pair $(G,P)$ in a uniform way.
These two terms correspond to the boundary terms of the functional $\theta_{G,L}$.
We further compute  a spherical Fourier coefficient, which
is used  to define the functional $\theta_{G,L}$.
In Section $6$ we formulate and prove our main Theorem \ref{main:thm}.
In Section $7$ we compute the functional $\theta_{D_3,D_2}$ for the group $SL_4$.
 In Section $8$ we compute the proportionality constant
between the  functionals $\theta_{G,L_1}$ and $\theta_{G,L_2}\circ T$
corresponding to different pairs $(G,L_1)$ and $(G,L_2)$, where $T$ is the transition map between two models.
This result is applied in Section $9$ to obtain  expansions
for automorphic functionals on the groups $E_6,E_7$ and $E_8$.

\section*{Acknowledgments}
We thank Dmitry Gourevitch for answering our questions and explaining the work \cite{Gourevitch}.
We  thank Hezi Halawi for his help in computation
of Gindikin-Karpelevich factors. We thank the referee for raising questions
regarding various results whose proofs in literature often
valid only for fields of characteristic zero, and for  sketching
the arguments to resolve them.

The first author is partially supported by ISF grant 1643/23.
The second author is partially supported by ERC grant 101142781.

\section{Preliminaries}\label{sec:preliminaries}
\subsection{The field $F$.}
Let  $k=\bF_q$ be a finite field, $C$ be a smooth, projective,
absolutely irreducible curve over  $k$ of genus $g\ge 1$ and  $F$ 
 be the field of  rational functions on the curve $C$.
 For any closed point $v\in C$ we denote
 by $F_v$ the local field with the ring of integers
$\mO_v$ and the residue field $k_v$.
We write $q_v = |k_v|$ and denote the norm by $|\cdot|_v$
and the valuation on $F_v$  by $val_v$. 
Then  $|\varpi_v|_v=q_v^{-1}$ for any generator
 $\varpi_v$ of the maximal ideal of $\mO _v$.

We fix an additive character $\psi:F\backslash \A$ that will
be used throughout the paper as follows. 
Let $\omega_C $ be  a non-zero regular differential form
on the curve $C$, that exists since the curve $C$ is of genus $\geq 1$. We denote 
by $\kappa _{ \omega_C} :\A \to k$ the homomorphism such that 
$$\kappa _{\omega_C}(f)=\sum_{v\in |C|} \Tr_v( 
\Res_v (f_v\omega_C)), \quad f=(f_v) \in \A$$
where $\Tr_v:k_v\rightarrow k$ is the trace form and $\Res_v$ is
the residue map with values in the residue field $k_v$.  
 
 We fix an additive character $\psi _0: k\to \C^\times $ and put
 $\psi := \psi _0\circ \kappa _{ \omega_C} $.
 It is clear that $\psi=\Pi_v \psi_v$
 defines a non-trivial character of $F\backslash \A$
 with trivial restrictions on $\mO _v$ for closed points $v$ of $C.$
 Let $a=(a_v)$ be the idele such that  for any $v$ the character
 $t_v\rightarrow \psi_v(a_v^{-1}t_v)$  has conductor $\mO_v$.
 One has $|a_v|_v\le 1$ for all $v$.

\sms

For any $v\in C$ we denote by $dx_v$ the Haar measure $dx_v$ on $F_v$ such that
$vol(\mO_v)=1$.
The measure  $dx=q^{1-g}\cdot\Pi_v dx_v$ is the Tamagawa measure on $\A$.

\sms

The absolute value  of the discriminant $\Delta_F$ of $F$ is $q^{2g-2}.$

 
\subsection{The divisor function}
Local divisor function $\sigma_v(m,s)$ for $m\in \Z$ and
$s\in \C$ is defined by  
$$\sigma_v(m,s)=\left\{\begin{array}{ll}\sum_{k=0}^{m} q_v^{ks} & m\ge 0 \\ 
0& m<0
\end{array}\right.$$
For any idele $b=(b_v)\in \A^\times $ define
$\sigma(b,s)=\prod_v\sigma_v(val_v(b_v),s),$
which is well-defined since $\sigma_v(0,s)=1$. 
In particular, for the idele $a$ associated to the character $\psi$
as above we have $\sigma(a,s)\neq 0$ for any integer $s$.
This fact is used in the proof of Proposition \ref{main:term:equality}.

\subsection{The completed zeta-function.}

The Hasse-Weil zeta function $\zeta(s)$ of $F$ has form 
$$\zeta(s)=\frac{\Pi_{i=1}^g (1-a_i q^{-s})(1-q a_i^{-1} q^{-s}) }{(1- q^{-s})(1- q^{1-s})},$$
where $|a_i|=q^{1/2}.$
This function is factorizable
$$\zeta(s)=\Pi_v \zeta_v(s), \quad \zeta_v(s)=\frac{1}{1-q^{-s}_v}$$ 
for all $v$.

The completed zeta-function, defined by  
$$\xi(s)=q^{-(1-g)s}\zeta(s)=q^{-(1-g)s}\cdot \Pi_v \zeta_v(s)$$ satisfies the functional equation 
$\xi(s)=\xi(1-s)$. It has simple poles at $s=1$ and $s=0$. 
We denote $\lim_{s\to 1} \xi(s)(s-1)=R$. 

\subsection{Group notation}
Let $G$ be a split  simply-connected Chevalley group,
 of type $A_n, D_n$ or $E_n$, defined over $\Z$.
We fix a  maximal split torus $T$ and a Borel subgroup
$B=T\cdot N$.
Let $R$ be the set of roots with respect to $T$.
The choice of the Borel subgroup defines the partition $R=R^+\cup R^-$
into positive and negative roots and the set
$\Delta=\{\alpha_i\}\subset R^+$ of simple roots.
We write $\alpha_0$ for the highest root. 
Similarly, we have  $R^\vee\supset \Delta^\vee$, the set of coroots.
The characters $\omega_\alpha, \alpha\in \Delta$ denote the
fundamental weights. We fix a Chevalley-Steinberg  system
$(e_\alpha:\G_a\rightarrow N_\alpha, \alpha\in R)$ of pinnings for $(G,T),$
where $N_\alpha$ is the root subgroup of $\alpha$.

The labeling for simple roots in the group $G$ is as follows: 
$$
 D_n \dynkin[labels={\alpha_1,,,\alpha_{n-2},\alpha_{n-1},\alpha_n},label directions={,,,right,,}]D{} \quad 
 E_6  \dynkin[labels={\alpha_1,\alpha_2,\alpha_3,\alpha_4,\alpha_5,\alpha_6}]E6\quad 
E_7 \dynkin[labels={\alpha_1,\alpha_2,\alpha_3,\alpha_4,\alpha_5,\alpha_6,\alpha_7}]E7\quad
E_8 \dynkin[labels={\alpha_1,\alpha_2,\alpha_3,\alpha_4,\alpha_5,\alpha_6,\alpha_7,\alpha_8}]E8 
$$
For any root $\alpha=\sum_{i=1}^m n_i \alpha_i\in R^+$
we denote $l(\alpha)=\sum_{i=1}^m n_i.$

Let $L\subset G$ be a  Levi subgroup in a maximal parabolic subgroup $P=LV$.
Then $G_1=[L,L]$ is simply-connected with a torus $T_1=T\cap [L,L]$.
We denote by $R_L$ the root system of $G_1$ with respect to $T_1$.
We also denote by $\beta_0$ the unique root in $\Delta-\Delta_L$.
There is a modular character
$2\rho_P=\sum_{\alpha\in R^+-R^+_L} \alpha$
that is trivial on $[L,L]$ and hence defines a character on $L$.
We write $|\rho_P|$ for $|2\rho_P|^{1/2}$.
The fundamental weight $\omega_{\beta_0}$ can be extended to a character of $L.$
We write $\Ind^{G}_P(s)$ for the normalized induction $\Ind^{G(F_v)}_{P(F_v)}|\omega_{\beta_0}|^s$.

We write $W, W_L$ for the Weyl groups of $G$ and $L$ respectively.
The group $W$ is generated by simple reflections $s_\alpha, \alpha\in \Delta$.
For any finite sequence $1\le k_i\le rk(G),$
the element $s_{\alpha_{k_1}}\ldots s_{\alpha_{k_r}}$ is denoted by $w(k_1,\ldots,k_r)$. The length of an element $w\in W$ is denoted by $l(w)$. The longest
element in $W$ and $W_L$  is denoted by $w_0$ and $w_0^L$ respectively. 

For any standard parabolic subgroup $P=LV$ we write $P^{op}$ for 
the opposite parabolic subgroup  with   Levi decomposition $P^{op}=LV^{op}$.

\subsection{The Eisenstein series on $SL_2$}
In this section let $G=SL_2$ with a Borel subgroup $B=TN$.  
The unique simple root is denoted by $\alpha$ and the non-trivial element 
of the Weyl group by $w$. 
The usual pinning gives identification $N(F)\simeq  F$ and
$T(F)\simeq  F^\times$.

Let $I_B(s)$ denote the representation $\Ind^G_B(|\rho_B|^s)$
of $SL_2(F_v),$ where $|\rho_B|^2$ is the modular character of the torus.
It is irreducible for $s\neq \pm 1$.

We recall basic facts about the Eisenstein series for $G=SL_2$.
For $f(\cdot,s)\in \Ind^{G(\A)}_{B(\A)}(|\rho_B|^s)$ the Eisenstein series is defined by 
$$\mE(f,s,g)=\sum_{\gamma\in  B(F)\backslash G(F)} f(\gamma g,s),$$
which converges for $Re(s)>1$. Let $f^0(g,s)$ be the spherical vector in $\Ind^{G(\A)}_{B(\A)}(|\rho_B|^s)$ 
satisfying $f^0(e,s)=1$. 

Let $\mM(s)\in \Hom_{G(\A)}
\left(\Ind^{G(\A)}_{B(\A)}(|\rho_B|^s), \Ind^{G(\A)}_{B(\A)}(|\rho_B|^{-s})\right)$ be the 
standard intertwining operator defined by 
$$\mM(s)(f)(g)=\int\limits_{N(\A)} f(wng,s) dn, $$
which converges for $Re(s)\gg 1$ and has meromorphic continuation to
the whole complex plane.  
Here the measure $dn$ on $N(\A)$  is the transfer of Tamagawa measure on $\A$ 
via the pinning $e_\alpha:\A\rightarrow N(\A).$

The additive character $\psi$ of $\A$, trivial on $F$ defines the additive character of $N(\A)$
trivial on $N(F)$ which we also denote by $\psi$. 

The formula for the constant term and the
$(N,\psi)$ Fourier coefficients for the
spherical Eisenstein series $\mE(f^0,g,s)$  are given below.

\begin{Prop}\label{prop:Eis:SL2}
\begin{enumerate}
\item
$\mE(f,s,g)^{N}=f(g,s)+\mM(s)(f)(g,-s),$
\item
  $\mM(s)(f^0)(e,-s)=c(s)f^0(e,-s)$ where $c(s)=\frac{\xi(s)}{\xi(s+1)}.$
\item $\mE(f^0,s,e)^{N,\psi}=q^{-s(1-g)}\cdot
  \frac{\sigma(a,-s)}{\xi(s+1)}\neq 0$ for $s>0$,
where the idele $a$ is associated to the character $\psi$.
\end{enumerate}
\end{Prop}

\begin{proof} The first two parts are standard computations. 
For the last one 
$$\int\limits_{N(F)\backslash N(\A)}\mE(f^0,s,n)\overline{\psi(n)}dn=
q^{1-g}\prod_v \int\limits_{N(F_v)} f_v^0(wn)\psi_v(n)dn=$$
$$q^{1-g}\frac{\sigma(a,-s)}{\zeta(s+1)}=
q^{-(1-g)s}\frac{\sigma(a,-s)}{\xi(s+1)}$$
\end{proof}
\subsubsection{The operator $\mM'$}

The function $c(s)$ in Proposition \ref{prop:Eis:SL2}
is holomorphic at $s=0$. Its Taylor series near $s=0$
is $$c(s)=-1+A s+\ldots $$
for some $A\in \C$.  In particular, $\mM(0)+Id\in Aut_{G(\A)}(I_B(0))$
annihilates the spherical vector $f^0$ and hence vanishes.
The  leading term of $\mM(s)+Id$ in Laurent expansion around $s=0$
is denoted  by $\mM'$.

Define a $T(F)N(\A)$ invariant functional $f\mapsto \mM'(f)(e)$.

\begin{Lem}\label{M:der} For $t=\alpha^\vee(r)$ one has
  $\mM'(t\cdot f^0)(e)=|r|(2\ln_q(|r|)+A)$
\end{Lem}

\begin{proof}
One has $f^0(\alpha^\vee(r),s)=|r|^{s+1}f^0(e,s)$

\begin{multline*}
|r|^{-1}\mM(s)(t\cdot f^0)(e)=|r|^{-1}\int\limits_{N(\A)} f^0(wn\alpha^\vee(r),s)dn=\\
|r|^{-s}\int\limits_{N(\A)} f^0(wn,s)dn=|r|^{-s}\frac{\xi(s)}{\xi(s+1)}.
\end{multline*}
Hence
$$|r|^{-1}(\mM(s)+Id)(tf^0)=|r|^{-s}\frac{\xi(s)}{\xi(s+1)}+|r|^s.$$
%
%
Taking the derivative with respect to $s$ at $s=0$ we obtain
 $$\mM'(\alpha^\vee(r)f^0)(e)=|r|(2\ln_q|r|+A).$$

\end{proof}

\section{Minimal representations}\label{sec:minimal}
In this section $F$  denotes a local non-archimedean field of arbitrary characteristic and the subscript
$v$ is omitted from all the data.
Let us  review known models of  minimal representation $\Pi$ of $G(F)$
and mention their useful properties.
For every model there is a maximal parabolic
subgroup $P=LV$ whose action in this model is explicit. We start
with the description of the pairs $(G,L)$ to be considered in this paper.

\subsection{Admissible and weakly admissible pairs}
\begin{Def}\label{admissible:def}
The pair $(G,L)$ is called {\bf admissible} if either
\begin{itemize}
\item The group $G$ is of type $E_i$ for $i=6,7,8$
and the unipotent radical $V$ of the maximal parabolic subgroup $P=LV$
is either abelian or Heisenberg group, or  
\item The group $G$ is of type $D_n$ for $n\ge 4$ and 
$L$ is of type $D_{n-1}.$ In this case the group  $V$ is abelian. 
\end{itemize}
\end{Def}

When $V$ is abelian (resp. Heisenberg)
the pair $(G,L)$ is called of Siegel type, or $S$-pair for short
(resp. of Heisenberg type or $H$-pair)

Given an  admissible pair $(G,L)$ we complete it to a triple $(G,L,M)$, 
where $M$ is a maximal Levi subgroup of $L$, as follows. Let $\beta_0$ 
be the simple root defining $L$. There exists unique simple root $\beta_1$
connected to $\beta_0$ in the Dynkin diagram, 
which defines a maximal parabolic subgroup $Q=MU$ in $L$. 
We write $G_1=[L,L]$ and $Q_1=Q\cap G_1=M_1 U$.  The pair $(G_1,M_1)$
 is not necessarily  admissible. For the triple $(G,L,M)$ the pairs $(G,L)$ and $(G_1,M_1)$
 are both called {\bf weakly admissible.}
 Note that any weakly admissible pair of $H$-type is  admissible.

\subsubsection{The $(s,d)$-parameters of weakly admissible pairs}

Let $(G,L)$ be an admissible pair and $\Pi, \Pi_1$ be minimal 
representations of $G(F)$ and $G_1(F)$ respectively.
The representation  $\Pi$ (resp. $\Pi_{1}$) is known,
see \cite{KazhdanSavin}, Theorem $4$  and \cite{Savin:invent},
Corollary $4.2.$ to be  unique irreducible subrepresentation 
of a degenerate principal series induced from $P$ (resp. $Q_1$).
In other words, there exist numbers $s_0,s_1>0$
such that 
$$\Hom_{G(F)}(\Pi,\Ind^G_P |\omega_{\beta_0}|^{-s_0})\neq 0, \quad  
\Hom_{G_1(F)}(\Pi_{1},\Ind^{G_1}_{Q_1} |\omega_{\beta_1}|^{-s_1})\neq 0.$$
Equivalently, for the Jacquet modules $(\Pi)_V$ and $(\Pi_{1})_U$ one has 
by Frobenius reciprocity,
$$\Hom_{L(F)}(\Pi_V,|\omega_{\beta_0}|^{d_0+1})\neq 0, \quad 
\Hom_{M_1(F)}({\Pi_{1}}_U,|\omega_{\beta_1}|^{d_1+1})\neq 0,$$
where 
\begin{equation}\label{s0:d0,s1:d1}
  -s_0+\la \rho_P,\beta_0^\vee\ra=d_0+1, \quad
  -s_1+\la \rho_Q,\beta_1^\vee\ra=d_1+1.
\end{equation}
The pairs $(s_0,d_0)$ and $(s_1,d_1)$ are called $(s,d)$-parameters for the 
weakly admissible pairs $(G,L)$ and $(G_1,M_1)$ respectively. 

The restriction map defines isomorphism $\Ind^{L}_{Q}|\omega_{\beta_1}|^{-s_1}
\simeq \Ind^{G_1}_{Q_1}|\omega_{\beta_1}|^{-s_1}$ of $G_1(F)$ representations. 
Hence the representation $\Pi_{1}$ can be extended to  $L(F)$
as the unique unramified subrepresentation of $\Ind^{L}_{Q}|\omega_{\beta_1}|^{-s_1}$. 

Here is the list of all $(G,L,M)$ triples
together with the parameters $(s_0,d_0)$ for the pair $(G,L)$ and $(s_1,d_1)$
for the pair $(G_1,M_1)$. 

 \begin{equation}\label{list}
\begin{array}{|c|c|c|c|c|c|c|}
  \hline
  G & L & M & (s_0,d_0) & (s_1,d_1)&
  \la\rho_P,\beta^\vee_0\ra& \la\rho_Q,\beta^\vee_1\ra\\  \hline
  D_n, & D_{n-1} & D_{n-2} & (1,n-3) & (1,n-2) & n-1 & n-2\\ \hline
  E_6& D_5 & A_4 & (3,2 )& (2,1)& 6 & 4\\ \hline
  E_7 & E_6 & D_5 & ( 5,3)& (3,2 )& 9 & 6\\ \hline
  E_6 & A_5 &  A_2\times A_2& (7/2,1) & (2,0)& 11/2 & 3\\ \hline
  E_7 & D_6 & A_5 & (11/2,2)  & (3,1) & 17/2 & 5 \\ \hline
  E_8 & E_7 & E_6&  (19/2,4) & (5,3) & 29/2 & 9 \\ \hline
 \end{array}
\end{equation}
 For the top triple we assume $n\ge 4$. The pairs $(G,L)$ in three
 last rows are $H$-pairs. 

For any weakly admissible pair $(G,L)$ we describe below the model
for the minimal representation with an explicit action of $P=LV$.
We treat $S$-pairs and $H$-pairs separately. 

\subsection{The model $\mS(X(F))$ for an $S$-pair $(G,L)$}
Let $(G,L)$ be a weakly admissible  $S$-pair, i.e. $V$ is abelian.
The  models with an explicit action of $P=LV$ 
have been systematically  studied in \cite{Savin:invent}, and then further in
\cite{SavinWoodbury}. The ground field in these
papers is assumed to be $p$-adic, but this is not used in the
quoted results in any essential way.


The model depends on a non-trivial additive character $\psi$ of $F$.
The pinning of $G$ defines a basis of the group $V^{op}(F)$,
as a $F$-vector space. Let $\|\cdot\|$ be the norm on $V^{op}(F)$ given by the
maximum of $v$-adic norms of coordinates. 

Any element $x\in V^{op}(F)$ defines 
the character $\psi_x$ of $V(F)$ by
$\psi_x(v)=\psi(\la v,x\ra),$ where $\la\cdot,\cdot\ra$
is the Killing form. 

The starting point to construct the model is to show that
the $V$-spectrum of the minimal representation $\Pi$ is ``small''.
Let us denote by  $\Pi_{V,\psi_x}$ the space of $(V,\psi_x)$ coinvaraints .

\begin{Prop}\label{S:uniqueness:local}
  Let $X\subset V^{op}$ be the minimal non-zero $L$-orbit.
  Then
$$\Pi_{V,\psi_x}=\left\{ \begin{array}{ll}\C & x\in X(F)\\
    0 & {\rm otherwise}\end{array}\right.  $$
\end{Prop} 
\begin{proof}
  The proof is essentially contained in \cite{Weissman}. 
  Let $Q=MH$ denote the standard Heisenberg parabolic subgroup.
  The key observation in \cite{Weissman} is that
  there is a decomposition $V=N_1\cdot A$, where
  the group $A=H\cap V$ is a maximal abelian subgroup of $H$
  and $N_1=M\cap V$ is  the unipotent radical of $M\cap P$.
  Let  $x_0=e_{-\alpha_0}(1)\in X(F).$
  For any $L$ orbit in $V^{op}$ there exists a representative
  $x$ such that $\psi_x|_{A}=\psi_{x_0}|_A$.
 For such $x$ one has  $\Pi_{A,\psi_x}\simeq (\rho_\psi)_{A,\psi_x}$.
   Since $\rho_\psi\simeq \ind^H_A \psi_x$ it is easy to see
   that $(\rho_{\psi})_{A,\psi_x}$ is one-dimensional
   and $N_1$ acts trivially on it.
  The restriction  of $\psi_x$ to $N_1$ is trivial if and only if $x\in X(F)$. 
 Hence $\Pi_{V,\psi_x}\simeq \C$ for $x\in X(F)$ and is zero otherwise. 
We thank the referee for the  sketch of this proof. 
\end{proof}
 
  We follow the terminology of \cite{Weissman}.
  The $L(F)$-orbits in $V^{op}(F)$ and consequently the characters
  $\psi_x$ of $V(F)$ are assigned ranks. The $V$-rank of a representation
  $\pi$ is the maximum of ranks of $\psi_x$, such that $\pi_{V,\psi_x}\neq 0$.
  The elements in the minimal non-zero
  orbit $X(F)$ have rank one, hence  the minimal representation $\Pi$ has
  $V$-rank one.

Using Proposition \ref{S:uniqueness:local} one realizes
$\Pi$ in a space of functions on $X(F)$ as follows. 

The group $P(F)$ acts on the space of smooth functions
$\mS^\infty(X(F))$  via 
\begin{equation}\label{P:action:S:pair}
\left\{\begin{array}{ll}
    g\cdot \phi(x)=|\omega_{\beta_0}(g)|^{d_0+1} \phi(g^{-1}xg) & g\in L(F)\\
    v\cdot \phi(x)=\psi_x(v) \cdot \phi(x) & v\in V(F)
\end{array}
\right.
\end{equation}

This action preserves the subspace $\mS_c(X(F))$
of smooth functions of compact support
and is unitary with respect to the standard inner product,
given by the unique, up to constant, $G_1(F)$ invariant measure on $X(F)$.

By \cite{Savin:invent}, Theorem $6.5$, see also \cite{SavinWoodbury},
Section $3$, there exists a $P(F)$-equivariant embedding
$\Pi\hookrightarrow \mS^\infty(X(F)),$
whose image is denoted by $\mS(X(F)).$ The fact that $\Pi$ is of
$V$-rank one is essential in the proof.

Note that the functional $l_x:\mS(X(F))\rightarrow \C$
defined by $l_x(\phi)=\phi(x)$ spans the one-dimensional
space $\Hom_{V(F)}(\Pi,\psi_x)$.

\subsubsection{The normalized spherical function}
The space $\mS(X(F))$ contains a one-dimensional space of spherical,
i.e. $K_0=G(\mO)$-invariant
functions.
\begin{Prop}(\cite{SavinWoodbury})
  Let us assume that $\psi$ is of conductor $\mO$.
  \label{spherical:formula:S}
  \begin{enumerate}
    \item 
      The function $\phi^0\in \mS(X(F))$ given by
      $$\phi^0(x)=\sigma(m,d_0), \quad \|x\|=q^{-m}$$
  is spherical, satisfying   $\phi^0(x)=1$ for $\|x\|=1$.
It is called the normalized spherical function.  
\item If $d_0>0$, which is always true for admissible pairs, one has   
\beq \label{spherical:function:S}
\phi^0(x)=\zeta(-d_0)+\zeta(d_0)\|x\|^{-d_0}, \quad \|x\|\le 1
\eeq
and zero otherwise.
  \end{enumerate}
  \end {Prop}

In the adelic setting we will be forced to consider character
$\psi$ of arbitrary conductor. For any $\psi$ denote the model
by $\mS_\psi(X(F))$ stressing the dependence of the action of $V(F)$
on the character $\psi$. The map
$$r_a:\mS_\psi(X(F))\rightarrow \mS_{\psi_a}(X(F)), \quad
r_a(\phi)(x)=\phi(ax),$$
where $\psi_a(x)=\psi(ax),$ defines a $G(\A)$ equivariant isomorphism
between the models. In particular $r_a(\phi^0)$ defines
the normalized spherical function in the model $\mS_{\psi_a}(X(F))$.

\begin{Cor} Let $\psi_v$ be a character of conductor $a\mO_v$. 
 The function $\phi^0\in \mS(X(F))$ given by
      $$\phi^0(x)=\sigma(m+val(a),d_0), \quad \|x\|=q^{-m}$$
  is spherical, satisfying   $\phi^0(x)=1$ for $\|x\|=|a|^{-1}$.

If $d_0>0$ then 
\beq 
\label{spherical:function:transition}
r_a(\phi^0)(x)=\phi^0(ax)=\zeta(-d_0)+\zeta(d_0)\|ax\|^{-d_0}
\eeq where $\|x\|=q^{-m}\le |a|^{-1}$ and zero otherwise. 
\end{Cor}

If the pair $(G,L)$ is admissible, then the pair
$(G_1,M_1)$ is weakly admissible of $S$-type. 
Hence the construction above is also applicable
to the pair $(G_1,M_1)$ giving rise to model $\mS(X_1(F))$
for the minimal representation $\Pi_{1}$ of $G_1(F)$. We denote by $\phi^0_{1}$
the normalized spherical function in $\mS(X_1(F))$. 

\subsubsection{Jacquet module $\Pi_V$}

 The space $\mS(X(F))$ contains the space
of functions of compact support, is contained in the space of functions
of bounded support and  fits in the $P(F)$ equivariant exact sequence:
$$0 \rightarrow \mS_c(X(F))\rightarrow \mS(X(F))
\rightarrow \Pi_V\rightarrow 0.$$ Thus $\Pi_V$
is naturally identified with the space of germs $[\mS(X(F))]_0$ at zero. 
Let us describe $\Pi_V$ as a representation of $L(F)$.  

\begin{Prop}\label{jacquet:S}
  Let $(G,L)$ be an admissible pair. There is an
  isomorphism of $L(F)$-modules 
  $$\Pi_V\simeq |\omega_{\beta_0}|^{d_0+1}\oplus
  |\omega_{\beta_0}|^{1+d_0(1-1/\kappa)}\Pi_{1},$$
  where $\kappa=\la \omega_{\beta_0},\omega_{\beta_0^\vee}\ra$.
\end{Prop}

\begin{proof}  
Consider the decomposition $\phi^0={\phi^0}'+{\phi^0}'',$
where ${\phi^0}',{\phi^0}''\in \mS^\infty(X(F))$ are defined by
$${\phi^0}''(x)=\zeta(d_0)\|x\|^{-d_0},\quad  {\phi^0}'=\phi^0-{\phi^0}''.$$
Hence ${\phi^0}'(x)=\zeta(-d_0)$  in a neighborhood of $0$
and ${\phi^0}''(x)$ is a homogeneous function.  
Since $\phi^0$ generates $\mS(X(F))$, it follows
by Iwasawa decomposition $G=P\cdot K_0$ and formulas 
(\ref{P:action:S:pair}) for $P$-action that any function $\phi\in \mS(X(F))$
has unique decomposition
$$\phi=\phi'+\phi'',\quad  \phi',\phi''\in
\mS^\infty(X(F)),$$ where
\begin{itemize}
  \item
    $\phi'$ is a constant in a neighborhood of  $0$,  denoted by $\phi'(0).$ 
    The map $\phi\mapsto \phi'(0)$ belongs to
    $\Hom_{P(F)}(\mS(X(F)),|\omega_{\beta_0}|^{d_0+1})$
    
  \item $\phi''$ is a homogeneous function of degree $-d_0$
    in a neighborhood of $0$.  Since $G_1(F)$ acts transitively on $X(F)$, and the stabilizer of the element 
$x_0=e_{-\beta_0}(1)$ is $\Ker \omega_{\beta_1}|_{Q_1}\subset Q_1(F)$, 
the space of homogeneous functions of degree $-d_0$ is identified 
with the space of sections $\Ind^{G_1}_{Q_1}(-s_1),$
where 
\begin{equation}\label{d0:s1}-s_1+\la \rho_Q,\beta_1^\vee\ra =d_0.
\end{equation}

The map $\phi\mapsto \phi''$, belongs to
$\Hom_{G_1(F)}(\mS(X(F))_V, \Ind^{G_1}_{Q_1}(-s_1)),$
and the image is the subspace generated by a non-zero spherical function.
For admissible pair $(G,L)$ the number  $s_1$ is the $s$-parameter for the pair $(G_1,M_1)$
and the image is exactly the minimal representation of $G_1$.
Let $\lambda^\vee$  be the coweight generating the connected component
of $Z(L)$, proportional to $\omega_{\beta_0^\vee}$. Applying formula
for the action of  $\lambda^\vee$ we have 
$$\lambda^\vee(r)\cdot (\phi^0)''(e_{-\beta_0}(1))=
|\omega_{\beta_0}(\lambda^\vee(r))|^{1+d_0(1-1/\kappa)} (\phi^0)''(e_{-\beta_0}(1)).$$

\end{itemize}

To sum up the map $\phi\rightarrow (\phi'(0),\phi'')$ defines the isomorphism 
$$\Pi_{V}\simeq |\omega_{\beta_0}|^{d_0+1}\oplus |\omega_{\beta_0}|^{1+d_0(1-1/\kappa)}\otimes
\Pi_1$$ as  representations of $L(F)$. 
\end{proof}

The decomposition  $\phi=\phi'+\phi''$ in the proof above allows us
to define the following important ingredients

\begin{Def}
  \be
\item The functional
  $p_{X}\in \Hom_{P(F)}(\mS(X(F)),|\omega_{\beta_0}|^{d_0+1})$ defined
  by  $p(\phi)=\phi'(0)$. In particular, $p_X(\phi^0)=\zeta(-d_0)$.
\item The map $\mB$ in the one-dimensional space 
  $\Hom_{P(F)}(\mS(X(F)),\mS(X_1(F))),$ called boundary map
  is defined by $\mB(\phi^0)=\phi^0_1.$
\ee
\end{Def}

\begin{Remark}  
The Proposition \ref{jacquet:S} has also been proven in \cite{Savin:invent}
using Iwahori-Hecke algebra computations.
\end{Remark}

\subsection{The model $\mS(X(F))$ for an $H$-pair $(G,L)$}
Let $G$ be of type $E_6,E_7$ or $E_8$ and
$P=LV$ be the standard Heisenberg parabolic subgroup. Note that the pair
$(G,L)$ is always admissible.
The additive complex character $\psi$ of $F$
can be viewed, using the pinning,  as a character of the one-dimensional
center $Z(F)$ of $V(F)$.

A polarization of the symplectic space $V/Z$ gives rise
to the abelian subgroups $V_1,V_2\subset V$ with $[V_1,V_2]=Z$. 
The group $V^+=V_2Z$ is a maximal abelian subgroup of $V$.
Let $T_0(F)$ be the  one-dimensional torus  generated by $\beta_0^\vee.$
Define $X=T_0\cdot V^+$.
The group  $X(F)=T_0(F)V_1(F)$ acts freely and transitively
on the set of  characters $\psi$ of $V^+(F)$
whose restriction to $Z(F)$ is non-trivial.
Choosing $\psi_e$ as a base point, for any $x\in X(F)$ we obtain
a character $\psi_x$ of $V^+(F)$.

In the pioneering paper \cite{KazhdanSavin} the authors construct
a special unitary representation $\hat\Pi$, by extending
the unitary completion of $\ind^P_{[P,P]}\rho_\psi$ to an irreducible
unitary representation of $G(F)$. The space of
$\rho_\psi=\ind^V_{V^+}\psi$ is naturally identified with $\mS_c(V_1(F))$.
Hence the space $\ind^P_{[P,P]}\rho_\psi$ can be identified
with the space of functions of compact support $\mS_c(X(F)),$
so that $\hat \Pi$ is realized on the space $L^2(X(F))$
of square-integrable functions. 
 The subspace of smooth vectors in $\hat \Pi$, denoted by
$\mS(X(F))$, define the irreducible smooth representation
$\Pi$. It is easy to see that the subspace
$\Pi(Z)=Span\{\Pi(z)\phi-\phi, z\in Z,\phi\in \mS(X(F))\}$ equals 
$\mS_c(X(F))$.
 Let us show that $\Pi$ is the minimal representation in the sense
 of Definition
\ref{min:def}.
  Indeed, there is an exact sequence
  $$0\rightarrow \Pi(Z) \rightarrow\Pi\rightarrow \Pi_Z\rightarrow 0$$
  Applying to it the exact functor of $(Z,\psi)$ coinvariants we obtain
  $$\Pi_{Z,\psi}\simeq \Pi(Z)_{Z,\psi}\simeq
  (\ind^P_{[P,P]}\rho_\psi)_{Z,\psi}\simeq \rho_\psi$$
as $[P,P]$ representations. 
Since  $\rho_\psi\simeq \ind^V_{V^+}\psi$ it follows by Frobenius reciprocity
that
\begin{equation}\label{H:uniqueness:local:V+}
  \Hom_{V^+(F)}(\Pi,\psi)\simeq \C
\end{equation}
for any character $\psi$, whose restriction to $Z(F)$ is non-trivial.

This gives rise to the $P$-equivariant exact sequence
$$0\rightarrow \mS_c(X(F))\rightarrow \mS(X(F))\rightarrow \Pi_Z\rightarrow 0$$.

\subsubsection{The model $\mS(Y(F))$ for $\Pi_Z$}
Our next goal is to construct a model for the representation $\Pi_Z$.
We closely follow the approach in \cite{KazhdanPolishchuk},
where $\cha(F)$ can be zero or positive. 

Let $Y$ be a minimal non-zero $L$ orbit in $(V/Z)^{op}$.
As in the case of $S$-pairs, the
elements of $y\in (V/Z)^{op}(F)$ parameterize characters $\psi_y$ of $(V/Z)(F)$
or equivalently characters of $V(F)$. 
The action of $P$ on $\mS^\infty(Y(F))$ is defined by formulas
in (\ref{P:action:S:pair}) and $Z(F)$ acts trivially. 

The starting point of constructing a model is the following proposition
\begin{Prop}\label{H:uniqueness:local:V}
$$\Pi_{V,\psi_y}=\left\{ \begin{array}{ll}\C & y\in Y(F)\\
    0 & {\rm otherwise}\end{array}\right.  $$
\end{Prop} 

\begin{proof} It is Proposition $3.2.2$ in \cite{KazhdanPolishchuk}.
  A different proof, which also does not use the assumption
  $\cha(F)=0$, appears in \cite{GanSavin} Proposition $11.5$.
\end{proof}  

Using this result, it is shown in
Theorem $6.1.1$ \cite{KazhdanPolishchuk}
that there is an embedding  of $\Pi_Z$ in the
space of locally constant sections of the line bundle over $Y(F)$
whose fiber over any $y\in Y(F)$ is $\Pi_{V,\psi_y}$.
For groups of type $E_i, i=6,7,8$,
Corollary $3.3.2$ defines an a $P(F)$-equivariant embedding of 
$\Pi_Z$ in the space of locally constant functions $\mS^\infty(Y(F))$.
Note that  for groups of type $D_n$ there is no such embedding. 
The image, denoted by
$\mS(Y(F)),$ fits in the exact $P(F)$-equivariant sequence
$$0 \rightarrow \mS_c(Y(F))\rightarrow \mS(Y(F))\rightarrow
\Pi_V\rightarrow 0.$$
Thus $\Pi_V$ is naturally identified with the space of germs $[\mS(Y(F))]_0$
at zero. 

Clearly $\Hom_{P}(\Pi, \Pi_Z)$ is one-dimensional.
We fix a map $\iota$ in the one-dimensional space
$\iota:\Hom_{P(F)}(\mS(X(F)),\mS(Y(F))$ 
as in \cite{KazhdanPolishchuk}, section $3$.
 
\subsubsection{The normalized spherical function}
Since $\Pi$ is unramified, it contains a one-dimensional
space of spherical functions. We shall need a formula for
$\iota(\phi^0)$ as a function on $Y(F)$, for a properly normalized $\phi^0$,
which follows from \cite{KazhdanPolishchuk}, Theorem $1.1.3$.
together with formula $(3.2)$. 

\begin{Prop}\label{spherical:formula:H}
  Let $\psi$ be a character of conductor $\mO$.
  Let $\phi^0$ be a spherical function normalized by the condition
  $\iota(\phi^0)(y)=1$ for $\|y\|=1$. Then one has 
\beq \label{spherical:function:H}
\iota(\phi^0)(y)=\left\{\begin{array}{ll}
  \zeta(-d_0)+\zeta(d_0)\|y\|^{-d_0}=\sigma(m,d_0) & \|y\|=q^{-m}\le 1\\
  0& {\rm otherwise}
\end{array}
  \right.
  \eeq
\end{Prop}

As in the case of $S$-pairs we extend the formula for characters
of arbitrary conductor. 

\begin{Cor} 
 Let $\psi$ be a character of conductor $a\mO$.
  Let $\phi^0$ be a spherical function normalized by the condition
  $\iota(\phi^0)(y)=1$ for $\|y\|=|a|^{-1}$. Then one has 
\beq \label{spherical:function:H}
\iota(\phi^0)(y)=\left\{\begin{array}{ll}
\zeta(-d_0)+\zeta(d_0)\|ay\|^{-d_0}=
\sigma(m+val(a),d_0) & \|y\|=q^{-m}\le |a|^{-1}\\
  0& {\rm otherwise}
\end{array}
  \right.
  \eeq
\end{Cor}

\subsubsection{Jacquet module $\Pi_V$}

Finally, the description of Jacquet module $\Pi_V$
is similar to the case of $S$-pair.
 Note that for any $H$-pair one has
$\la\omega_{\beta_0},\omega_{\beta_0^\vee}\ra =2$. 

\begin{Prop} \label{jacquet:H}
  One has $\Pi_V\simeq |\omega_{\beta_0}|^{d_0+1}\oplus
  |\omega_{\beta_0}|^{1+d_0/2}\otimes \Pi_1$
  as representations of $L(F).$    
\end{Prop}

\begin{proof}
  The proof is the same as for Proposition \ref{jacquet:S}.
  We use the formula for  normalized spherical function
  and Iwasawa decompostion for $G$ to write for  any
  function $\phi\in \mS(X(F))$
  \beq \label{jacquet:decomp:H}
  \iota(\phi)=\iota(\phi)'+\iota(\phi)'',
  \eeq
  where $\iota(\phi)'\in \mS^\infty(Y(F))$ is constant near zero and
  $\iota(\phi)''\in \mS^\infty(Y(F))$ is homogeneous near zero.
  Each summand gives rise to the summand of $\Pi_V$.
  See also  \cite{KazhdanPolishchuk}, Theorem
  $6.2.2$ or \cite{Savin:invent} Theorem $4.1$. 
\end{proof}

Since  $(G_1,M_1)$ is a weakly admissible $S$-pair,
the representation $\Pi_1$ has a model in a space $\mS(X_1(F))$
defined earlier with the normalized spherical function $\phi^0_1$.

Similar to $S$-pair case we define the functional $p_Y$
and the boundary map $\mB$ as follows. 

\begin{Def}
  \be
  \item
We define the functional $p_Y$ on $\mS(Y(F))$ by
$p_Y(\iota(\phi))=\iota(\phi)'(0)$,
where $\iota(\phi)'$ is taken from the decomposition (\ref{jacquet:decomp:H}).
\item
The map $\mB$ in the one-dimensional space  
$\Hom_{P(F)}(\mS(X(F)),\mS(X_1(F)))$ is fixed by
$\mB(\phi^0)=\phi^0_1.$
\ee
\end{Def}

\subsection{Relation between parameters}
The following relations between the parameters 
will be used for the proof of main theorem. 

\begin{Cor}\label{parameter:relations}
  For a triple $(G,L,M)$ from the table \ref{list}
  with the parameters $(s_0,d_0)$ and  $(s_1,d_1)$ holds
  $$s_0-s_1=\la \rho_P,\beta_0^\vee\ra-\la \rho_Q,\beta_1^\vee\ra-1,
  \quad  d_0=d_1+1.$$
\end{Cor}

\begin{proof} The identities  follow immediately from the relations
$$-s_0+\la \rho_P,\beta_0^\vee\ra=d_0+1,\quad -s_1+\la \rho_Q,\beta_1^\vee\ra=d_1+1$$
$$d_0=-s_1+\la \rho_Q,\beta_1^\vee\ra $$
mentioned above
\end{proof}

\subsection{The transition map}
Let $(G,L_1)$ and $(G,L_2)$ be weakly admissible pairs.
Denote the models of the minimal representation $\Pi$ of $G(F)$ by  
$\mS(X_1(F))$, and  $\mS(X_2(F))$ 
and the normalized spherical functions by $\phi^0_1, \phi^0_2$ respectively. 
The transition map $T\in \Hom_{G(F)}(\mS(X_1(F)),\mS(X_2(F))$
is uniquely defined by the condition $T(\phi^0_1)=\phi^0_2.$

\section{The global model}\label{sec:global}
In this section $F$ is a global field and $\A$ is its ring of adeles.
Let $\psi=\Pi_v\psi_v:F\backslash \A\rightarrow \C^\times$ be the fixed complex
additive character  whose associated idele $(a_v)$
satisfies $|a_v|_v\le 1$ for all places $v$. 
Let $(G,L)$ be a weakly admissible pair with parameters $(s_0,d_0)$. 
We define the global model, associated to the pair $(G,L)$ and the character
$\psi$, for the representation
$\Pi=\otimes_v \Pi_v$ of $G(\A)$ as the restricted tensor product
$\mS(X(\A))=\otimes'_v\mS(X(F_v))$ with respect to the
normalized spherical functions $\phi^0_v$ above.
We write $\phi^0=\otimes \phi^0_v$, the global spherical function.

For an   $S$-pair  $(G,L)$ let us fix
 the element $x_0=e_{-{\beta_0}}(1)\in X(F)$ with $\|x_0\|=1,$
 Combining Propositions \ref{spherical:function:S},
 and  equation
 (\ref{spherical:function:transition}) we get
 $$\phi^0(x_0)=\sigma(a,d_0)\neq 0.$$
 
 For an   $H$-pair $(G,L)$ let us fix
 the element  $y_0=e_{-{\beta_0}}(1)\in Y(F)$ with $\|y_0\|=1$).
 Combining Proposition  \ref{spherical:function:H} and equation
 (\ref{spherical:function:transition})
 we get
$$\iota(\phi^0)(y_0)=\sigma(a,d_0)\neq 0.$$ 
 
\subsection{The functionals  $\theta_X$, $\theta_Y$} 

For the variety $X,Y$ the functionals $\theta_X$ on $\mS(X(\A))$  and $\theta_Y$
on $\mS(Y(\A))$ are defined by
$$\theta_X(\phi)=\sum_{x\in X(F)}\phi(x), \quad \theta_Y(\phi)=\sum_{y\in Y(F)}\phi(y).$$
For any $\phi$ the sums are finite, and hence the functionals are well-defined.

\subsection{The regularized functionals $p_X,p_Y$}
Let $(G,L)$ be an admissible  $S$-pair.
The map $p_X\in \Hom_{P(\A)}(\mS(X(\A)), |\omega_{\beta_0}|^{d_0+1})$
associates to every function its regularized value at $0$.
Precisely,  let $\phi=\otimes_{v}\phi_v\in \mS(X(\A)$
and $S$ be a finite set of places such that $\phi_v=\phi^0_v$ for $v\notin S$. 
  \begin{equation}\label{p:X:global}
  p_X(\phi)=\prod_{v\in S} p_{X,v}(\phi_v)\cdot \zeta^{S}(-d_0),
 \end{equation}
  where $\zeta^S(s)$ is the partial zeta-function.
  Since $p_{X,v}(\phi^0_v)=\zeta_v(-d_0)$,
  the functional $p_X$ is well-defined. 

  Similarly, for $H$-pairs,
  the map $p_Y\in \Hom_{P(\A)}(\mS(Y(\A)), |\omega_{\beta_0}|^{d_0+1})$
  is defined by
\beq
p_Y(\iota(\phi))=\prod_{v\in S} p_{Y,v}(\iota_v(\phi_v))\cdot \zeta^{S}(-d_0).
\eeq

\subsection{The boundary map}
The global boundary map $\mB=\otimes \mB_v$ is defined by
\beq \label{B:global}
\mB\in \Hom_{P(\A)}(\mS(X(\A)),\mS(X_1(\A)),\quad \mB(\phi^0)=\phi^0_1.
\eeq

\section{Eisenstein series}\label{sec:eisenstein}

Fix a weakly admissible pair $(G,L)$. The representation
$\Pi$ of $G(\A)$ is unique irreducible quotient of the degenerate principal
series $\Ind^G_{P'}(s_0),$ where  $P'=L'V',$ is the standard parabolic subgroup
with the Levi factor  $L'=w_0(L)$.
The degenerate Eisenstein series $\mE_{P'}(f,g,s)$ admits
a simple pole at $s=s_0$ and the residues span the automorphic
representation $\bar\theta(\Pi)$ isomorphic to $\Pi$,
contained in the space of square-integrable automorphic functions
$\mA_2(G)$. The functional $\varphi\mapsto \varphi(e)$
is $G(F)$-invariant. The goal of this section is to decompose
this functional as sum of three (resp. four) terms for $S$-pairs
(resp. $H$-pairs) and to determine the normalizing constant $c_{G,L}$
which is used in definition of the automorphic functional
$\theta_{G,L}$ on $\mS(X(\A))$.

Recall the notation
$$\{\beta_0\}=\Delta-\Delta_{L}, \quad\{\beta_0'\}=\Delta-\Delta_{L'},\quad
\Ind^G_{P'}(s)=\Ind^{G(\A)}_{P'(\A)}|\omega_{\beta_0'}|^{s}.$$
The induced representations contains a spherical section $f^0(g,s)$ normalized by $f^0(e,s)=1$.
For any flat section $f(\cdot,s)$ define the Eisenstein series by
$$\mE_{P'}(f,s,g)=\sum_{\gamma \in P'(F)\backslash G(F)}f(\gamma g,s),$$
which converges for $Re(s)>\la\rho_P,\beta_0^\vee\ra$ and
has a meromorphic continuation for the complex plane.

\subsection{Poles of $\mE_{P'}(f,s,g)$}
At  $s=\la\rho_P,\beta_0^\vee\ra$ the Eisenstein series $\mE_{P'}(f,s,g)$
has at most  simple pole. The pole is attained by $f=f^0$ and
the residue is a constant function on $G(F)\backslash G(\A)$.
The point $s_0$ is the second largest singular point. We summarize
the results regarding the residual representation below. 

\begin{Thm} 
  \be
  \item
    Let $(G,L)\neq (A_5,A_2\times A_2)$ be a weakly admissible pair that appears in the list
    \ref{list} either as
    $(G,L)$ or $(L,M).$ 
  The Eisenstein series $\mE_{P'}(f,s,g)$ admits at most simple pole at $s=s_0$.
The pole is attained by the spherical function. 
\item
  The function $\bar\theta_{G,L}(f)(g)=Res_{s=s_0}\mE_{P'}(f,s,g)$
  is square-integrable. The functions $\bar\theta_{G,L}(f)$, as
  $f$ runs over standard sections in $\Ind^{G}_{P'}(s)$
  span an irreducible automorphic representation $\bar\theta(\Pi)$  
of $G(\A)$, isomorphic to $\Pi$.
\item The multiplicity of $\Pi$ in the discrete square-integrable
  automorphic spectrum is one. 
\ee    
\end{Thm}

\begin{proof}
The poles of degenerate Eisenstein series are governed
by the poles of its constant term along the unipotent radical
of a maximal parabolic subgroup, which can be explicit
computed for the spherical function. Eventually, this has been done
for every pair separately. 

The proof for the pairs $(E_8,E_7)$, $(E_7,E_6),$ $(E_6,D_5)$ and
$(D_n,A_{n-1})$ appears in \cite{GRS}.
The proof for $(E_7,D_6)$ appears in \cite{Gurevich:thesis}.
The proof for $(E_6,A_5)$ can be deduced from M.Sc. Thesis of H. Halawi,
\cite{halawi2022poles}
as well as for other pairs with $G=E_6$ or $E_7$. 
The case $(D_n,D_{n-1})$ with $n\ge 3$ can be done similarly. 

The proof for part $(3)$ for groups of type $D_n$ and $E_n$ appears 
\cite{KazhdanPolishchuk}, Corollary $8.1.3$. 

\end{proof}

For the remaining pair $(G,L)=(A_5,A_2\times A_2)$
the Eisensentein series has a simple pole,
and the residual representation is irreducible, isomorphic to
the minimal representation, as shown in \cite{HanzerMuic}.
The only difference is that the resulting representation is
not contained in the space of square-integrable automorphic functions.

\subsection{Fourier expansion of $\bar\theta_{G,L}(f)$}
Let  $(G,L)$ be a pair of  $S$-type.
The Fourier expansion of the function
$\bar\theta_{G,L}(f)$ along
the abelian group $V(F)\backslash V(\A)$
reads
\begin{equation}\label{Fourier:S}
\bar \theta_{G,L}(f)(e)=\sum_{x\in X(F)} \bar\theta_{G,L}(f)^{V,\psi_x}(e)+
\bar \theta_{G,L}(f)^V(e).
\end{equation}
For $(G,L)$ of $H$-type we consider first the Fourier expansion along
$Z$. The non-trivial Fourier coefficients are further expanded along
the group $V_+$. The trivial Fourier coefficient is expanded along
 the abelian group $V/Z$. 
\begin{equation}\label{Fourier:H}
  \bar\theta_{G,L}(f)(g)=\sum_{x\in X(F)}
  \bar\theta_{G,L}(f)^{V^+,\psi_x}(g)+
  \sum_{y\in Y(F)}\bar \theta_{G,L}(f)^{V,\psi_y}(g)+
  \bar\theta_{G,L}(f)^V(g).
\end{equation}

It follows from Propositions
\ref{S:uniqueness:local}, \ref{H:uniqueness:local:V} that other Fourier coefficients vanish.

\subsection{The constant term $\bar\theta_{G,L}(\cdot)^V$}\label{sec:CT}
In this subsection we assume that $(G,L)$ is an {\sl admissible pair}
and describe the residue  of the constant term $Res_{s=s_0}\mE_{P'}(f,s,g)^V$.

By a standard computation the constant term along $V$ equals:
\begin{equation}\label{const:term:along:V:notation}
  \mE_{P'}(f,s,g)^V=\sum_{w\in \Psi_{L',L}}  J_w(f,s,g), \quad J_w(f,s,g)=\mE^L_{Q_w}(\mM_w(s)f,s_w,g),
  \end{equation}
where
\begin{enumerate}
\item $\Psi_{L',L}=\{w\in W: w(\Delta_L), w^{-1}(\Delta_{L'})\subset R^+\}$
  is the set of shortest representatives in $W_{L'}\backslash W/W_L.$
\item $Q_w=w^{-1}P' w\cap L$  is either $L$ or its maximal
  parabolic subgroup.
\item $\chi_{P',s}=|\omega_{\beta_0'}|^s\cdot |\rho_{P'}|^{-1}\cdot |\rho_B|$
  so that $\Ind^G_{P'}(s)\hookrightarrow \Ind^G_B\chi_{P',s}.$
\item $s_w$ is such that $w^{-1}(\chi_{P',s})|_L=\chi_{Q_w,s_w}$
\item $$\mM_{w}(s):\Ind^G_{P'}(s)\rightarrow \Ind^L_{Q_w}(s_w), 
\quad \mM_w(s)(f)(g)=\int\limits_{(N\cap w^{-1}Nw)(\A)\backslash N(\A)}f(w ng,s) dn$$
  are the standard intertwining operators, restricted to $L$. 
\item $\mM_{w}(s)(f^0)=c_w(s)f_1^0,$ where
$$c_w(s)=\Pi_{\alpha\in R_w}
  \frac{\xi(\la\chi_{P',s},\alpha^\vee\ra)}
       {\xi(\la\chi_{P',s},\alpha^\vee\ra+1)},$$
where $R_w=\{\alpha\in R^+, w^{-1}(\alpha)\in R^-\}$
       and $f_1^0$ is the normalized spherical section in
       $\Ind^L_{Q_w}(s_w)$.
\end{enumerate}

\subsubsection{The distinguished elements $v_0,v_1$}
The constant term  $\bar\theta_{G,L}(f)^V$ is always a sum of two terms
which eventually give rise to two boundary terms in Theorem \ref{main:thm}.
This subsection is devoted to study these terms and their properties. 

We begin with the standard parabolic subgroups $P=LV\subset G$ defined by the root
$\beta_0$ and $Q=MU\subset L$ defined by the root $\beta_1$. Let $w_0,w_0^L,w_0^M$
be the longest  elements in Weyl groups of $G,L,M$ respectively.

Let $P'=L'V'\subset G$ be the standard parabolic,  with $L'=w_0(L)$
defined by the root $\beta_0'=-w_0(\beta_0)$, and 
$Q'=M'U'\subset L$ be the standard parabolic with $M'=w_0^L(M),$
defined by the root $\beta'_1=-w_0^L(\beta_1)$.

Let $\alpha_0$ be the highest root in $R$ and $\alpha^L_0$
be the highest root in $R_L$. 
It is easy to check the following Lemma
\begin{Lem} Let $\tilde\beta=s_{\beta_0}(\alpha_0),
  \tilde \beta_L=s_{\beta_1}(\alpha_0^L)$.
  \be
\item If $(G,L)$ is an $S$-pair then $\tilde\beta=\alpha_0$.
\item If $(G,L)$ is an $H$-pair then $\tilde\beta=\alpha_0-\beta_0$.
\item $\tilde\beta_L=\alpha_0^L$.
  \item\label{length:tilde:beta}
 $ w_0(\tilde\beta)=-\tilde\beta, \quad l(\tilde\beta)=2\la\rho_P,\beta^\vee_0\ra-1,
  \quad
  l(\tilde\beta_L)=2\la\rho_Q,\beta^\vee_1\ra-1.$
 
    \ee
\end{Lem}

Let us define the distinguished elements in $W$
$$  v_0=w_0 w_0^L, \quad v_0^L=w_0^L w_0^M, \quad v_1=v_0 s_{\beta_0}(v_0^L)^{-1}.$$

\begin{Lem}\label{v0:v1}
\begin{enumerate}
\item $v_0\in \Psi_{L',L}$ and $v_0^{-1}P'v_0\cap L=L,$
\item $v_0(\beta_0)=-\tilde\beta,$
\item $v_1\in \Psi_{L',L}$ and $v_1^{-1}P'v_1\cap L=Q',$
\item $v_1^{-1}(\chi_{P',s_0})|_{G_1}=\chi_{Q_1',s_1}.$
\end{enumerate}
\end{Lem}

\begin{proof}
\begin{enumerate}
\item The first par is easy.
\item 
 It suffices to show that $w_0^L(\beta_0)=\tilde\beta$.
The elements $w^L_0(\beta_0),\tilde \beta$ both have form
$\beta_0+\sum_{\alpha\in \Delta_L} n_\alpha\cdot  \alpha$
for some non-negative integers $n_\alpha$.
Hence $w^L_0(\beta_0)-\tilde\beta\in Span(\Delta_L)$. Let us
show that $\la w^L_0(\beta_0),\alpha^\vee\ra=\la\tilde\beta,\alpha^\vee\ra$
for all $\alpha\in \Delta_L$ and hence $w^L_0(\beta_0)=\tilde\beta.$

One has $\Delta_L=\Delta_{M'}\cup \beta_1'$. For any $\alpha\in \Delta_L$
one has
$$\la w^L_0(\beta_0),\alpha^\vee\ra=\la \beta_0,w_0^L(\alpha^\vee)\ra=\left\{
\begin{array}{ll} 1 & \alpha=\beta_1'\\
0 & {\rm otherwise}\end{array}\right.$$

Assume $(G,L)$ is an $S$-pair.
The root $\beta_1'=-w_0^L(\beta_1)$ is the simple root, corresponding
to the Heisenberg parabolic subgroup. Hence
$\la\tilde\beta,\alpha^\vee\ra=\la \alpha_0,\alpha^\vee\ra=1$
for $\alpha=\beta_1'$ and zero otherwise.

If $(G,L)$ is $H$-pair then $\beta_1'=\beta_1$.
Hence $\la\tilde\beta,{\beta_1'}^\vee\ra=\la\alpha_0-\beta_0,\beta_1^\vee\ra=1$
and  $\la\tilde\beta,\alpha\ra=0$ for $\alpha\in \Delta_{M'}=\Delta_M$. 
Thus  $w^L_0(\beta_0)=\tilde\beta$, as required.
\item 
One has  $\Delta_L=\Delta_{M'}\cup \beta_1'$ 
and $v_1=w_0w_0^L s_{\beta_0} w_0^L w_0^{M'}$.
We check directly
$$v_1(\Delta_{M'})\subset \Delta_{L'}, \quad 
v_1(\beta'_1)=\tilde\beta-\tilde\beta_{L'}\in R^+-R^+_{L'}$$
Hence $v_1^{-1}P'v_1\cap L=Q'$.

To prove $v_1\in \Psi_{L',L}$ we need to verify that
$v_1(\Delta_{L'}), v_1^{-1}(\Delta_{L'})\subset R^+$.
The first containment follows from above. 
We write $v_1^{-1}=w_0^{M'} w_0^Ls_{\beta_0}w_0^Lw_0$. 
One has $w_0^Lw_0 (\Delta_{L'})=\Delta_L=\Delta_M\cup \{\beta_1\}$.
Further
$$w_0^{M'} w_0^Ls_{\beta_0}(\Delta_M)=\Delta_{M'},\quad  w_0^{M'} w_0^Ls_{\beta_0}(\beta_1)=
\tilde\beta-\tilde\beta_{L'}\in R^+,$$
as required.

 \item 
We already know that $v_1^{-1}(\chi_{P',s_0})|_{G_1}=\chi_{Q_1',t}$ for some $t$. 
To show that $t=s_1$ we compute the pairing with $\beta_1'^\vee$.
$$\la v_1^{-1}(\chi_{P',s_0}),\beta_1'^\vee\ra=
\la \chi_{P',s_0},\tilde\beta^\vee-\tilde\beta^\vee_{L'}\ra=
s_0-\la\rho_P,\beta_0\ra+2\la\rho_Q,\beta^\vee_1\ra,$$
since   $l(\tilde\beta-\tilde\beta_{L'})=
2(\la\rho_P,\beta^\vee_0\ra-\la\rho_{Q'},{\beta'}_1^\vee\ra)$
by (\ref{length:tilde:beta}).

On the other hand  $\la\chi_{Q_1',s_1},\beta_1'^\vee\ra=
s_1+\la\rho_{Q_1'},{\beta'}_1^\vee\ra -1.$
Hence 
$$v_1^{-1}(\chi_{P',s_0})|_{G_1}=\chi_{Q_1',s_1}\Leftrightarrow
s_0-s_1=\la\rho_P,\beta^\vee_0\ra-\la\rho_{Q'},{\beta'}^\vee_1\ra-1$$
which is true by Corollary \ref{parameter:relations}.
\end{enumerate}
\end{proof}

\begin{Thm}\label{const:term}
One has for any standard section $f\in\Ind^G_{P'}(s)$ 
$$\bar\theta_{G,L}(f)^V=\bar\mM_{v_0}(s)(f)(e)+\bar\theta_{G_1,M_1}(\mM_{v_1}(s_0)(f)),$$
where 
$\bar \mM_{v_0}=Res_{s=s_0}\mM_{v_0}(s)$. 
\end{Thm}

\begin{proof}

For the triples $(G,L,M)$ in the table \ref{list}  we find out by case-by-case consideration
that the terms $J_w(f,s,g)$ are holomorphic at $s=s_0$ unless
$w\in \{v_0,v_1\}$. Let us study the terms $J_{v_0}$ and $J_{v_1}$ in more detail. 

The operator $\mM_{v_0}(s)$ has a simple pole at $s=s_0$. Since $\Ind^G_{P'}(s_0)$
is generated by the spherical vector, it is enough to show that the factor
$c_{v_0}(s)$ has a simple pole at $s=s_0$. This can be seen 
from the following table. 

\begin{equation}
\begin{array}{|c|c|c|c|}
  \hline
  (G,L) & c_{v_0}(s)& s_0 & d_0\\
\hline  
(D_n,D_{n-1})& \frac{\xi(s)\xi(s-n+2)}{\xi(s+1)\xi(s+n-1)}
& 1 & n-3  \\
\hline
(E_6,D_5)& \frac{\xi(s - 2)\xi(s - 5)}{\xi(s + 6)\xi(s + 3)}
 &  3 &2 \\
\hline
(E_6,A_5)&
\frac{\xi(s - \frac{3}{2})\xi(s - \frac{5}{2})\xi(s - \frac{9}{2})\xi(2 s)}
 {\xi(s + \frac{5}{2})\xi(s + \frac{7}{2})\xi(s + \frac{11}{2})\xi(2 s + 1)}
  & 7/2 & 1\\
\hline
(E_7,E_6)&
\frac{\xi(s)\xi(s - 4)\xi(s - 8)}{\xi(s + 1)\xi(s + 5)\xi(s + 9)}
  & 5  & 3\\
\hline
(E_7,D_6)&
\frac{\xi(s - \frac{5}{2})\xi(s - \frac{9}{2})\xi(s - \frac{15}{2})\xi(2 s)}
     {\xi(s + \frac{7}{2})\xi(s + \frac{11}{2})\xi(s + \frac{17}{2})\xi(2s+1)} & 11/2 & 2\\
\hline
(E_8,E_7)& \frac{\xi(s - \frac{9}{2})\xi(s - \frac{17}{2})\xi(s - \frac{27}{2})\xi(2s)}
    {\xi(s + \frac{11}{2})\xi(s + \frac{19}{2})\xi(s + \frac{29}{2})\xi(2+1)}
 & 19/2 & 4\\
\hline
\end{array}
\end{equation}

 The operator $\mM_{v_1}(s)$ is holomorphic at $s=s_0.$
This follows from the holomorphicity of the factor $c_{v_1}(s)$ at $s_0$.
The values appear in the table below
$$
\begin{array}{|c|c|c|}
  \hline
  (G,L) & v_1 & c_{v_1}(s_0)\\
  \hline
  (D_n,D_{n-1}) & w(1) & \frac{\xi(n-1)}{\xi(n)}  \\
  \hline
  (E_6,D_5)&  w(65431) & \frac{\xi(4)}{\xi(9)} \\
  \hline
  (E_6,A_5) & w(24315436542) &
  \frac{\xi(3) \xi(4)\xi(5)}{\xi(6)\xi(8)\xi(9)}\\
  \hline
  (E_7,E_6)&  w(7654234567) &\frac{\xi(5)\xi(9)}{\xi(10)\xi(14)} \\
  \hline
  (E_7,D_6)& w(13425436542765431)& 
\frac{\xi(4)\xi(6)\xi(8)}
     {\xi(9)\xi(12)\xi(14)} \\
     \hline
 (E_8,E_7) & w(8765432143546 257 6453412345678) &
  \frac{\xi(6)\xi(10)\xi(14)}{\xi(15)\xi(20)\xi(24)}\\
       \hline
\end{array}
$$

 Hence the operator
 $\mM_{v_1}(s_0)(f):\Ind^G_{P'}(s_0)\rightarrow \Ind^L_{Q'}(s_1)$
 is well-defined. 
 
The Eisenstein series $\mE^L_{Q'}(f,s,e)$ has a simple pole at $s=s_1$
 and $$\bar\theta_{G_1,M_1}(f)(e)=Res_{s=s_1}\mE^L_{Q'}(f,s,e)|_{G_1(\A)}$$

It follows from Lemma \ref{v0:v1} that 
$$\bar\theta_{G,L}(f)^V(e)=\bar\mM_{v_0}(s_0)(f)(e)+
\bar\theta_{G_1,M_1}(\mM_{v_1}(s_0)f)$$
as required. 
\end{proof}

 \begin{Cor} For any weakly admissible pair the function
      $\frac{c_{v_0}(s)}{\xi(s-s_0-d_0)}$ admits a simple pole at $s=s_0$.
 \end{Cor}

\begin{proof}
  For admissible pairs $c_{v_0}(s)$ has a simple pole at $s=s_0$
  and $d_0>0$, hence $\xi(s-s_0-d_0)$ is holomorphic at $s=s_0$. 
Let us add the table for the values of $c_{v_0}(s)$ for weakly-admissible,
but not admissible pairs. 
For the pairs $(A_5,A_2\times A_2)$ and $(D_3,D_2)$ the factor
$c_{v_0}(s)$ has a pole of order $2$ at $s=s_0$.
These are exactly the cases, where
$d_0=0$  and  function $\frac{c_{v_0}(s)}{\xi(s-s_0-d_0)}$
still has a simple pole at $s=s_0$.

\begin{equation}
\begin{array}{|c|c|c|c|}
  \hline
  (G,L) & c_{v_0}(s) & s_0 & d_0\\
\hline
(A_5,A_2\times A_2)&
\frac{\xi(s)\xi(s-1)\xi(s-2)}{\xi(s+1)\xi(s+2)\xi(s+3)}
    & 2 & 0\\ 
\hline
(D_3,D_{2}) &
\frac{\xi(s - 1)\xi(s)}{\xi(s + 1)\xi(s + 2)}
  & 1 & 0\\
\hline
(D_5,A_{4}) &
\frac{\xi(s - 1)\xi(s - 3)}{\xi(s + 2)\xi(s + 4)}
 & 2 & 1 \\
\hline
(D_6,A_5) &
\frac{\xi(s - 2)\xi(s)\xi(s- 4)}{\xi(s + 1)\xi(s + 5)\xi(s + 3)}
 & 3 & 1 \\
\hline
 \end{array}
\end{equation}
 \end{proof}
 
\subsection{A  Fourier coefficient of the spherical function}
Let $(G,L)$ be an $S$-pair.
In order to define an automorphic functional $\theta_{G,L}$
we need to fix  an isomorphism $\tau_X$ in the
one-dimensional space $\Hom_{G(\A)}(\mS(X(\A)),\bar\theta(\Pi))$,
such that $\tau_X(\phi)^{V,x}=\phi(x)$ for all $x\in X(F)$.
It is enough to fix $\tau_X$ on one vector. 
The functions $\tau_X(\phi^0)$ and  $\bar \theta_{G,L}(f^0)$ are both spherical
 in $\bar\theta(\Pi)$ and hence are proportional. We write 
$\bar\theta_{G,L}(f^0)=c_{G,L}\cdot \tau_X(\phi^0)$, and expect
the equality
$$\bar\theta(f^0)^{V,\psi_{x_0}}(e)=
c_{G,L}\cdot \tau_X(\phi^0)^{V,\psi_{x_0}}=c_{G,L}\cdot \sigma(a,d_0)\neq 0$$
for the base point $x_0=e_{-\beta_0}(1)\in X(F)$.
The goal of this subsection is to determine the constant $c_{G,L},$
by computing the spherical Fourier coefficient $\bar\theta_{G,L}(f^0)^{V,\psi_{x_0}}(e)$.
This is the most involved computation in the paper. Similarly
for $H$-pair we need to compute   $\bar\theta_{G,L}(f^0)^{V,\psi_{y_0}}(e)$
for the base element $y_0=e_{-\beta_0}(1)\in Y(F)$.

\begin{Thm}\label{spherical:Fourier}
  Let $c_{G,L}=q^{(1-g)d_0}\cdot Res_{s=s_0} \frac{c_{v_0}(s)}{\xi(s-s_0-d_0)}.$ 
  \be
    \item For $S$-pair
\begin{equation}\label{S:c:GL:def}
\bar\theta_{G,L}(f^0)^{V,\psi_{x_0}}(e)=c_{G,L}\cdot \sigma(a,d_0).
\end{equation}

\item For $H$-pair
 \begin{equation}\label{H:c:GL:def}
\bar\theta_{G,L}(f^0)^{V,\psi_{y_0}}(e)=c_{G,L}\cdot \sigma(a,d_0).
\end{equation}
\ee
  \end{Thm}

\begin{proof}
  The proof, that treats $S$-pairs and $H$-pairs simultaneously,
  occupies the rest of the section. 
In the proof we repeatedly use the root $\beta_0$, and 
do not use  $\beta_1$. To ease notation we write $\beta$ for $\beta_0$. 
  
The character $\psi_{x_0}$ of $V(F)\backslash V(\A)$ can be
extended to a character of $N(F)\backslash N(\A),$ 
where $N$ is the unipotent radical  of
 the Borel subgroup $B$, by defining it to be trivial on all the root
 subgroups $R_\gamma$ for $\gamma\in R^+_L$. 
 The extended character is  also denoted by $\psi_{x_0}$.
 Our main tool is Theorem   $A$ in  \cite{Gourevitch}.
$$\bar\theta_{G,L}(f^0)^{V,\psi_{x_0}}=\bar\theta_{G,L}(f^0)^{N,\psi_{x_0}}.$$
Let $L_{\beta}$ be the Levi subgroup of $L$ generated by
the roots $\pm \beta$. Its Borel subgroup is denoted by
$B_{\beta}=T N_{\beta}$.

Let $N'\subset N$ be the subgroup
generated by all the positive root subgroups, except for $\beta$.
Then $N=N'\cdot N_{\beta}$. Clearly, 
$$\mE_{P'}(f^0,s)^{N,\psi_{x_0}}=\left(\mE_{P'}(f^0,s)^{N'}\right)^{N_\beta,\psi_{x_0}}.$$ 
We compute first  the constant term  $\mE_{P'}(f^0,s)^{N'},$
that is an automorphic function on $L_\beta.$
    
      $$\mE_{P'}(f^0,g,s)^{N'}=
      \sum_{v\in \Psi_{L',L_\beta}} J_v(s,g), \quad J_v(s,g)=
      c_v(s)\cdot \mE^{L_\beta}_{v^{-1}L'v\cap L_\beta}(f^0,s_v,g).$$

The term $Res_{s=s_0}J_v(s,g)$ can possibly have a non-zero $(N_\beta,\psi_{x_0})$
coefficient only if $c_{v}(s)$ has a pole at $s=s_0$ and  $v^{-1}L'v\cap L_\beta=B_\beta.$
The last condition is equivalent to $v(\beta)\in R^+-R^+_{L'}$.
Such Weyl element $v$
is called {\sl relevant}.
We show below that $v_0s_\beta$ is  {\sl unique} relevant element 
in $\Psi_{L',L_\beta}$.

\begin{Lem} The element $v_0s_\beta$ is relevant
  and $$Res_{s=s_0} J_{v_0s_\beta}^{N_\beta,\psi_{x_0}}(s)=q^{(1-g)d_0} \cdot
  Res_{s=s_0}\frac{c_{v_0}(s)}{\xi(s-s_0-d_0)}\cdot 
  \sigma(a,d_0).$$
\end{Lem}

  \begin{proof}    
Clearly $v_0 s_\beta(\beta)=\tilde\beta\in R^+-R^+_{L'}$.     

      It follows from (\ref{s-s0-d0}) that
      $$\la (v_0s_\beta)^{-1}(\chi_{P,s}),\beta^\vee \ra=
      \la \chi_{P,s},\tilde \beta\ra=s-s_0-d_0$$
     Hence
     $$J_{v_0s_\beta}(s)=c_{v_0s_\beta}(s)\cdot\mE^{L_\beta}_{B_\beta}(f^0,e,s-s_0-d_0).$$

      One has $[L_\beta,L_\beta]\simeq SL_2$. By Proposition \ref{prop:Eis:SL2}, part $3$, one has 
      $$\mE^{L_\beta}_{B_\beta}(f^0,s)^{N,\psi}=q^{-(1-g)s}\cdot
      \frac{\sigma(a,-s)}{\xi(s+1)}$$
  
     
      Hence 
      $$J_{s_{v_0s_\beta}}(s)^{N_\beta,\psi_{x_0}}= 
      q^{-(1-g)(s-s_0-d_0)}\cdot c_{v_0s_\beta}(s)\cdot 
      \frac{\sigma(a,d_0+s_0-s)}{\xi(s-s_0-d_0+1)}$$
      On the other hand $R_{v_0}=R_{v_0s_\beta}\cup \{\tilde \beta\},$
      hence
      $$c_{v_0}(s)=c_{v_0s_\beta}(s)\cdot
      \frac{\xi(s-s_0-d_0)}{\xi(s-s_0-d_0+1)}.$$
  So  $Res_{s=s_0}J_{v_0s_\beta}(s)^{N,\psi_{x_0}}$ equals   
      $$
      q^{(1-g)d_0}Res_{s=s_0}\frac{c_{v_0s_\beta}(s)}{\xi(s-s_0-d_0+1)}\cdot \sigma(a,d_0)=
      q^{(1-g)d_0}Res_{s=s_0}\frac{c_{v_0}(s)}{\xi(s-s_0-d_0)}\cdot \sigma(a,d_0),$$        
as required.
  \end{proof}

  It remains to show that there are no other relevant terms in $\Psi_{L',L_\beta}$.
  If $(G,L)$ is one of the pairs $(D_n,D_{n-1}), (A_5,A_2\times A_2),(D_5,A_4),(D_6,A_5)$
  then for all $v\neq v_0s_\beta\in \Psi_{L',L_\beta}$
    one has $v(\beta)\in R^+_L$  and hence $v$ is not relevant. 

For remaining pairs that involve the group of type $E$ 
the set $\Psi_{L',L_\beta}$ can be very large.
Hence we offer an inductive argument for these pairs,
that are all admissible, which still requires  case-by-case investigation. 

Let us describe the main idea first and then make the detailed computation
for $(G,L)=(E_6,D_5)$.

For any $(G,L)$ we choose a weakly admissible pair 
$(\hat L, \hat M)$ 
where $\hat L \subset G$ is a maximal Levi subgroup, different from $L$
and  $\{\beta\}=\Delta_{\hat L}-\Delta_{\hat M}$. 
This gives standard parabolic subgroups $\hat P=\hat L\hat V$ and 
$\hat Q=\hat M\hat U\subset \hat L.$
We claim that there exists element $u\in \Psi_{L',\hat L}$ such that 
$c_u(s)$ is holomorphic at $s=s_0$ and 
\begin{equation}\label{relevant:unique}
\bar\theta_{G,L}(f^0)^{N,\psi_{x_0}}=
c_{u}(s_0)\bar\theta_{\hat L,\hat M}(f_1^0)^{N_{\hat L},\psi_{x_0}},
\end{equation}
where $N_{\hat L}$ is  the unipotent radical of the
  Borel subgroup $B_{\hat L}$.

In particular, if the pair $(\hat L, \hat M)$ has unique relevant element
$v$ then  $uv$ is the unique relevant element for the pair $(G,L)$.

To prove (\ref{relevant:unique}) we compute the degenerate Fourier
coefficient in stages. Since
$N=\hat V\cdot N_{\hat L}$,
$$\bar\theta_{G,L}(f^0)^{N,\psi_{x_0}}=(\bar\theta_{G,L}(f^0)^{\hat V})^{N_{\hat L},\psi_{x_0}}$$

$$\bar \theta_{G,L}(f^0)^{\hat V}(g)=
\sum_{w\in \Psi_{L', \hat L}} I_w(f^0,g),\quad 
I_w(f^0,g)=Res_{s=s_0} \mE^{\hat L}_{Q_w}(\mM_{w}(s)f^0,s_w,g)$$
where $Q_w=w^{-1}P'w\cap \hat L.$
The sum  typically contains only two non-zero terms $I_{u_1}$ and $I_u$. However
$I_{u_1}(f^0,g)$ is a constant function and hence its $(N_{\hat L},\psi_{x_0})$
coefficient is zero. The factor $c_u(s)$ is holomorphic at $s=s_0$
and $u^{-1} P'u\cap \hat L=\hat Q$, 
where $s_w$ in the $s$-parameter for the pair $(\hat L,\hat M)$. Hence 
$I_u(f^0,g)=c_u(s_0)\cdot \bar\theta_{\hat L,\hat M}(f^0)(g),$
and (\ref{relevant:unique}) follows. 

Let us make detailed computation for $(G,L)=(E_6,D_5).$ The other cases are similar. 
We use the labeling of the roots in the Dynkin diagrams as in section \ref{sec:preliminaries}.
We denote by $P_i=L_iV_i$ the maximal parabolic subgroup of $G$
such that $\{\alpha_i\}=\Delta-\Delta_{L_i}$. The labeling of roots
  in $L_i$ is inherited from $\Delta$.  
  The maximal parabolic $Q_{i,j}=M_{i,j}U_{i,j}$ in $L_i$
  satisfies $\{\alpha_j\}= \Delta_{L_i}-\Delta_{M_{i,j}}$. For any
  Levi subgroup $L$ we denote by $N_L$ the unipotent radical of its
  Borel subgroup $B_L$.

One has $s_0=3$. We fix $L=L_6$ so that  $L'=L_1$ and fix $\hat L=L_1.$ 
Consequently $\hat M=M_{1,6}.$
We compute the constant term $\mE_{P_1}(f^0,g,s)^{\hat V}$ where $\hat V=V_1$.
     One has  $$\Psi_{L_1,L_1}=\{e,u_1=w(1),u=w(13425431)\}$$
     and it is clear that  $I_e(f^0,g)=0$. 
Let us show that $I_{u_1}(f^0)^{N_{\hat L},\psi_{x_0}}$ is zero. 
The factor $c_{u_1}(s)=\frac{\xi(s+5)}{\xi(s+6)}$ is holomorphic at $s=3$
and $$u_1^{-1}P_1u_1\cap L_1=Q_{1,3} \quad u_1^{-1}(\chi_{P',s})=\chi_{Q_{1,3},s+1}.$$
For $s=3$ one has $s+1=4=\la\rho_{Q_{1,3}}^{L_1},\alpha_3^\vee\ra$, 
i.e. $I_{u_1}(f^0,s,\cdot)$ is a constant function. 
Consequently, its non-trivial Fourier coefficient is zero. 

The factor $c_{u}(s)=\frac{\xi(s+2)\xi(s-1)}{\xi(s+3)\xi(s+6)}$
is also holomorphic at $s=3$. One has 
$$u^{-1}P_1u\cap L_1=Q_{1,6}, \quad u^{-1}(\chi_{P,s})=\chi_{Q_{1,6},s-2}$$
 For $s=3$ one has $s-2=1$, that is the $s$-parameter for the pair $(L_1,L_{1,6})$
 of type $(D_5,D_4)$, for which the uniqueness of the relevant element
 was already established. 

We tabulate for every pair $(G,L)$ the choice of $\hat L$ and the element $u$ in the table below. 
$$
\begin{array}{|c|c|c|}
\hline
  (G,L)   &    (\hat L, \hat M) & u\in \Psi_{L',\hat L}\\
  \hline
  (E_6,D_5)   & (D_5,D_4) &  w(13425431) \\ 
\hline
(E_6,A_5) & (D_5,A_4) & w(24354265431)\\
\hline
(E_7,E_6) & (D_6,D_5) &
w(76542314354265431)\\
\hline
(E_7,D_6) & (E_6,D_5) &
w(13425431654234567)\\
\hline
(E_8, E_7) & (D_7,D_6) & w(87654231435426543 1 76542 3456 87 65423143542 65431)
\\
\hline
\end{array}
$$

This finishes the proof of Theorem \ref{spherical:Fourier}.

\end{proof}

\section{The functional $\theta_{G,L}$ and Main Theorem}\label{sec:main:thm}
Let $(G,L)$ be a weakly admissible pair.
The spaces $\mS(X(\A))$ and $\bar\theta(\Pi)$ are isomorphic as
$G(\A)$ representations. We fix a $G(\A)$-equivariant map
$$r_X\in \Hom_{G(\A)}(\Ind^G_{P'}(s_0),\mS(X(\A))), \quad
r_X(f^0)=c_{G,L}\cdot \phi^0,$$
 and the isomorphism of $G(\A)$-representations 
$$\tau_X\in \Hom_{G(\A)}(\mS(X(\A)),\bar\theta(\Pi)), \quad
\tau_X (c_{G,L}\cdot  \phi^0)=\bar\theta_{G,L}(f^0),$$
where the constant $c_{G,L}$ is taken from Theorem \ref{spherical:Fourier}.
Obviously $\tau_X\circ r_X=\bar \theta_{G,L}$.

\begin{Def} For any weakly admissible pair $(G,L)$
  the automorphic functional $\theta_{G,L}$ of $\mS(X(\A))$
  is defined by $\theta_{G,L}(\phi)=\tau_X(\phi)(e)$.
  Equivalently, $\theta_{G,L}(r_X(f))=\bar\theta_{G,L}(f)(e)$. 
\end{Def}

In particular, for an admissible pair $(G,L)$,  completed
to the triple $(G,L,M),$ the automorphic functional
$\theta_{G_1,M_1}$ is defined on $\mS(X_1(\A))$ which
is the model for the minimal representation $\Pi_1$ of $G_1(\A)$.

We arrive to our main Theorem

\begin{Thm}\label{main:thm} Let $(G,L)$ be an admissible
  pair completed to the triple $(G,L,M)$ from table \ref{list}. 
  \be
  \item If $(G,L)$ is $S$-pair then 
\beq \label{main:S}
\theta_{G,L}=\theta_X+p_X+|\Delta_F|^{1/2}\theta_{G_1,M_1}\circ \mB.
\eeq
\item 
  If $(G,L)$ is $H$-pair then
\beq \label{main:H}
  \theta_{G,L}=|\Delta|^{-\frac{n+3}{2}}
  \theta_X+\theta_Y\circ \iota +p_Y\circ \iota+
  |\Delta_F|^{1/2}\theta_{G_1,M_1}\circ \mB.
  \eeq
\ee
\end{Thm}

The proof of Theorem occupies the rest of this section. It
follows from comparing terms in (\ref{main:S}) and (\ref{main:H})
with terms in the Fourier expansion
of  $\tau_X(\phi)$ in (\ref{Fourier:S}), (\ref{Fourier:H}).

We start with the main term.
\begin{Prop}\label{main:term:equality}
  \begin{enumerate}
  \item
 If $(G,L)$ is an  $S$-pair, then      
\beq \label{fourier:coeff:value:S}
 \tau_X(\phi)^{V,\psi_x}(e)=\phi(x),
 \quad \forall x\in X(F),\phi\in \mS(X(\A))
 \eeq
\item 
  If $(G,L)$ is  an $H$-pair, then  
\beq  \label{fourier:coeff:value:H1}
  \tau_X(\phi)^{V,\psi_y}(e)=\iota(\phi)(y),
  \quad \forall y\in Y(F),\phi\in \mS(X(\A))
  \eeq 
and there exists a constant $a_{G,L}$ such that
\beq \label{fourier:coeff:value:H2}
\tau_X(\phi)^{V_+,\psi_x}(e)=a_{G,L} \cdot \phi(x),
\quad \forall x\in X(F),\phi\in \mS(X(\A)).
\eeq

  \end{enumerate}
\end{Prop}

\begin{proof} For $S$-pair $(G,L)$
the functionals $\phi\mapsto \tau_X(\phi)^{V,\psi_{x}}(e)$
  and $\phi\mapsto \phi(x)$ belong to the one-dimensional space
  $\Hom_{V(\A)}(\Pi, \psi_{x})$ and hence are proportional by
  Proposition \ref{S:uniqueness:local}.
  It remains to show that there exists a vector on which that are equal
 and non-zero.  Recall that $\phi^0(x_0)=\sigma(a,d_0)\neq 0$.
  
 The group $L(F)$ acts transitively on $X(F)$(\cite{SavinKobayashi},
 Theorem $4.2$ and Corollary $4.3$).
For any $x\in X(F)$ there exists $g\in L(F)$ such that $g^{-1}xg=x_0$.
\begin{multline}
\tau_X(g\phi^0)^{V,\psi_x}(e)=
c_{G,L}^{-1} \bar\theta_{G,L}(f^0)^{V,\psi_{x}}(g)=\\
c_{G,L}^{-1}\cdot \bar\theta_{G,L}(f^0)^{V,\psi_{x_0}}(e)=
\sigma(a,d_0)=\phi^0(x_0)=(g\phi^0)(x).
\end{multline}
and hence equation (\ref{fourier:coeff:value:S}) holds for all
$\phi\in \mS(X(\A)), x\in X(F)$.

The proof of equation (\ref{fourier:coeff:value:H1}) for $H$-pairs is the same.
 It relies on the fact that $L(F)$ acts transitively on $Y(F)$ that is
 proved in \cite{KazhdanPolishchuk}, Lemma $8.1$
 and on Proposition \ref{H:uniqueness:local:V}.
 
 Let us prove the equality (\ref{fourier:coeff:value:H2}).
 The  functionals $\phi\mapsto \tau_X(\phi)^{V_+,\psi_x}(e)$ and
 $\phi\mapsto \phi(x)$
 are proportional by (\ref{H:uniqueness:local:V+}).
 Let $a_{G,L}$ be the constant of proportionality for $x_0=e$.
 For any $x\in X(F)$ by change of variables it is easy to see that
 for any automorphic form  $\varphi$ one has 
 $\varphi^{V_+,\psi_x}(x^{-1})=\varphi^{V_+,\psi_{x_0}}(e).$ Hence
 for any $\phi\in \mS(X(\A))$ and  $x\in X(F)$ one has
\begin{multline}
\tau_X(\phi)^{V^+,\psi_x}(e)=
 \tau_X(x\phi)^{V^+,\psi_x}(x^{-1})=\tau_X(x\phi)^{V,\psi_{x_0}}(e)=\\
 a_{G,L} \cdot (x\phi)(x_0)=a_{G,L}\cdot \phi(x).
 \end{multline}
Hence $a_{G,L}$ is the  proportionality constant for any $x$.  
\end{proof}

\begin{Cor}
For $S$-pairs we have
  \begin{equation}
\theta_X(\phi)=\sum_{x\in X(F)} \tau_X(\phi)^{V,\psi_x}(e)
     \end{equation} 

For $H$-pairs we have
  \beq
  a_{G,L}\cdot\theta_X(\phi)=\sum_{x\in X(F)}\tau_X(\phi)^{V_+,\psi_x}(e), \quad
  \theta_Y(\iota(\phi))=\sum_{y\in Y(F)} \tau_X(\phi)^{V,\psi_y}(e)
  \eeq

\end{Cor}

Next we compare the boundary terms in (\ref{main:S}) and (\ref{main:H}) with
the constant term  $\tau_X(\cdot)^V(e).$

\begin{Prop}
  \be
\item
  \beq
  \bar\mM_{v_0}(s_0)=\left\{ \begin{array}{ll} p_X\circ r_X & S-pair\\
    p_Y\circ \iota\circ r_X & H-pair
  \end{array}\right.
  \eeq

\item For pairs of both types one has 
  $$\theta_{G_1,M_1}\circ \mB\circ r_X=
  \bar\theta_{G_1,M_1}\circ \mM_{v_1}(s_0)(\cdot)(e)$$
  \ee
\end{Prop}

\begin{proof}
\be
\item   
  Let $(G,L)$ be an admissible  $S$-pair. In particular $d_0>0$
  and $\zeta(-d_0)$ is well-defined. Let us show that both
  functionals agree on $f^0$. Indeed, 
\begin{multline}
p_X\circ r_X(f^0)=c_{G,L}\cdot \zeta(-d_0)=
q^{(1-g)d_0}
\frac{\bar c_{v_0}(s_0)}{\xi(-d_0)}\cdot \zeta(-d_0)=\\
\bar c_{v_0}(s_0)=\bar\mM_{v_0}(s_0)(f^0).
\end{multline}

The functionals $p_X\circ r_X$ and $\bar\mM_{v_0}(s_0)$
belong to $\Hom_{P(\A)}(\Ind^{G}_{P'}(s_0),|\omega_{\beta_0}|^{d_0+1})$
and factor through $r_X$, hence proportional, by Proposition
\ref{S:uniqueness:local}.
Since they agree on $f^0$ and are non-zero, they are equal. 

The proof  that $p_Y\circ \iota\circ r_X= \bar\mM_{v_0}(s_0)$
for $H$-pairs is the same. 

\item  Both maps factor through $r_X$ and belong
  to $\Hom_{V(\A) G_1(F)}(\Ind^G_{P'}(s_0), \Pi_1)$ and hence
  are proportional. It is enough to show that they agree on $f^0$.

By definition 
  $$\theta_{G_1,M_1}\circ \mB\circ r_X(f^0)=
  |\Delta_F|^{1/2} c_{G,L}\cdot \theta_{G_1,M_1}(\phi^0_1)(e)$$
  and 
  $$\bar\theta_{G_1,M_1}\circ \mM_{v_1}(s_0)(f^0)(e)=
  c_{G_1,M_1}\cdot c_{v_1}(s_0)\theta_{G_1,M_1}(\phi^0_1)(e)$$

  Let us prove that
  \beq \label{second:identity}
  |\Delta_F|^{1/2} c_{G,L}= c_{G_1,M_1}\cdot c_{v_1}(s_0)\eeq

  Recall that   $v_1\cdot v_0^L=v_0 s_{\beta_0}$
and $R_{v_0}=R_{v_0s_{\beta_0}}\cup \{\tilde\beta\}$. One has 
  \begin{equation}\label{s-s0-d0}
  \la \chi_{P',s},\tilde \beta\ra =
  s+\la\rho_P,\beta_0^\vee\ra-l(\tilde\beta)=s-\la\rho_P,\beta_0^\vee\ra-1=
  s-s_0-d_0
  \end{equation}
by identities (\ref{length:tilde:beta})
 and  (\ref{s0:d0,s1:d1}).
$$c_{v_0}(s)=c_{v_0s_\beta}(s)\cdot \frac{\xi(s-s_0-d_0)}{\xi(s-s_0-d_0+1)}$$
Besides, $d_1=d_0-1$ by Corollary \ref{parameter:relations}. 
Hence one has
  $$c_{G,L}= q^{(1-g)d_0}Res_{s=s_0}\cdot \frac{c_{v_0}(s)}{\xi(s-s_0-d_0)}=
 q^{(1-g)d_0}\cdot Res_{s=s_0} \frac{c_{v_0s_\beta}(s)}{\xi(s-s_0-d_0+1)}=$$
  $$ q^{1-g}\cdot c_{v_1}(s_0) q^{(1-g)d_1} Res_{s=s_0}
  \frac{c_{v_0^L}(s-s_0+s_1)}{\xi(s-s_0-d_1)}=$$
  $$q^{(1-g)} \cdot  c_{v_1}(s_0)\cdot q^{(1-g)d_1}
  Res_{s=s_1}\frac{c_{v_0^L}(s)}{\xi(s-s_1-d_1)}=
  |\Delta_F|^{-\frac{1}{2}}\cdot c_{v_1}(s_0)\cdot c_{G_1,M_1},$$
  as required.

To summarize,  for $S$-pairs one has  
$$\theta_{G,L}=\theta_X+p_X+|\Delta_F|^{\frac{1}{2}}\theta_{G_1,M_1}\circ \mB$$
and for $H$-pairs one has
$$\theta_{G,L}=a_{G,L}\theta_X+\theta_Y\circ \iota+p_Y\circ \iota+
|\Delta_F|^{\frac{1}{2}}\theta_{G_1,M_1}\circ \mB$$

It follows from  \cite{KazhdanPolishchuk}, Theorem $8.2.1$ that
$a_{G,L}=|\Delta_F|^{-\frac{n+3}{2}}$. 
\ee
\end{proof}

\section{The pair $(D_3,D_2)$}\label{sec:D3:D2}
For the group $G_n$ of type $D_n,$ $n\ge 4$ we denote by $X_n$
the cone of non-zero  isotropic vectors in the split
quadratic space of dimension $2n-2$. The minimal representation
$\Pi_n$ of $G_n(\A)$ affords a model $\mS(X_n(\A)$.
Theorem \ref{main:thm} can be applied several times for the
functional $\theta_{D_n,D_{n-1}}$ to conclude that it equals:

$$\sum^{n-4}_{k=0}  |\Delta_F|^{\frac{k}{2}} \theta_{X_{n-k}}\circ \mB_{n}^{(k)}+
\sum_{k=0}^{n-4} |\Delta_F|^{\frac{k}{2}} p_{n-k}\circ \mB_n^{(k)}+ 
|\Delta_F|^{\frac{n-3}{2}}\theta_{D_3,D_2}\circ \mB_n^{(n-3)},$$

with the last term in the expression is $\theta_{D_3,D_2}$, where $\Pi_3$ is
a minimal representation of the group $D_3$, defined below,
although in this case it is not unique.   

The pair  $(D_3,D_2)$ is  weakly admissible, but not admissible, so
Theorem \ref{main:thm} does not hold in this case.
There are two related phenomena  that change the situation. 
\begin{itemize} 
\item Locally, the Jacquet module of the minimal representation $(\Pi_{3,v})_V$
is no longer a direct sum, but rather a non-split extension of a representation $\Pi_{2,v}$
of $L(F_v)$ by a one-dimensional character. The representation $\Pi_{2,v}$ is reducible of length $3$
and so is not the minimal representation of $G_1(F_v).$
\item Globally, the representation $\Pi_3$ is still realized automorphically as a residue of the Eisenstein series
$\mE_P(f,g,s),$ which has a simple pole the $s=1$.  The constant term along $V$ is a sum of two terms that 
are non-holomorphic at $s=1$, both having a pole of order $2$ and only their sum has a simple pole. 
\end{itemize}

We shall formulate and prove an analogue of Theorem \ref{main:thm}
for the pair $(D_3,D_2)$.
Through  rather long preparation 
that is needed to formulate it, we state several claims, all concern representations of groups of type $D_n$ with $n\le 3$. 
The claims are either fairly simple or standard and the proofs are left to  reader.

First let us introduce some notation. 
 The simply connected group $G$ of type $D_3$ is isomorphic to $SL_4$.
 The Levi subgroup $L$ of type $D_2$ is 
$$ L\simeq \{(g_1,g_2)\in GL_2\times GL_2: \det(g_1)\cdot \det(g_2)=1\}.$$
The simple roots and fundamental weights of $G$ are denoted by
$\alpha_i, \omega_i$ for $1\le i\le 3$.
We denote by $B=N\cdot T$ the Borel subgroup of $G$ and by $B_L=N_L\cdot T$ the Borel subgroup of $L$. 
The maximal split torus $T$ is  product of the subgroups $T_i$,
where $T_i$  is the image of $i$-th coroot. 
Any character of $T(F_v)$ or $T(\A)$  of the form
$\Pi^3_{i=1} |\omega_i|^{z_i}$ is denoted by $\chi_{(z_1,z_2,z_3)}$. 
The character $\omega_2$ is extended to the character $|\det|$ of $L$,
i.e. $\omega_2(g_1,g_2)=|\det(g_1)|$.

Recall that  $X_n$ is  the cone of isotropic non-zero vectors
in the $2n-2$-dimensional split quadratic space.
In particular $X_2$ is the variety of non-zero isotropic vectors in a two-dimensional
split quadratic space $(V_2,q_2)$ . Fixing a basis of $V_2$ so that 
$V_2(F)\simeq F\oplus F$ and $q_2(x,y)=xy$ we write  
$$X_2(F)=\{(b,0),(0,b), \quad  b \in F^\times \}\subset F\oplus F.$$

The representation $\Pi_3$ is defined to be
unique irreducible quotient of the degenerate principal series
$\Ind^G_P (1)$. This representation is minimal. 
In particular, it affords the model $\mS(X_3(\A))$
with the normalized spherical function  $\phi^0_3$.

Let us define several maps: 
\begin{itemize}
\item $\bar\theta_{G,L}(f)(g)=Res_{s=1}\mE_P(f,g,s)$
    \item The isomorphism
      $\tau_3\in \Hom_{G(\A)}(\mS(X_3(\A),\bar \theta(\Pi_3))$
      such that
      $$\tau_3(\phi)^{V,\psi_x}=\phi(x), \quad \phi\in \mS(X_3(\A)), x\in X_3(F).$$
    \item   $r_3:\Ind^G_P(1)\rightarrow \mS(X_3(\A))$ satisfying
      $\bar r_3(f^0)=\frac{R}{\xi(2)\xi(3)}\phi^0_3.$
    \item The functional $\theta_{G,L}$ on $\mS(X_3(\A))$ is defined by
      $\theta_{G,L}(\phi)=\tau_3(\phi)(e).$
\end{itemize}
In particular the following diagram is commutative
$$
\begin{diagram}
\node{\Ind^{G}_{P}(1)}
\arrow{e,t}{\bar\theta_{G,L}} 
\arrow{s,l}{r_3}
  \node{\bar\theta(\Pi_3)}
  \arrow{s,r}{\varphi\mapsto \varphi(e)}
  \\
\node{\mS(X_3(\A))}
\arrow{ne,t}{\tau_3}
\arrow{e,t}{\theta_{G,L}}
\node{\C}
\end{diagram}
$$

Next we shall define two families of non-zero Fourier coefficients of $\varphi\in \bar\theta(\Pi_3).$

For  any $x\in X_3(F)$ we have defined the character $\psi_x$ of $V(F)\backslash V(\A).$
For any $y\in X_2(F)$ we associate a character $\psi_y$ of $N(F)\backslash N(\A)$ as follows. 
For $y=(b,0)$ the character  $\psi_y$ is trivial on all the root subgroups except $N_{\alpha_1}$
and $\psi_y(e_{\alpha_1}(r))=\psi(yr).$ Similarly, for $y=(0,b)$ the character  $\psi_y$ is trivial on all the root subgroups except $N_{\alpha_3}$ and $\psi_y(e_{\alpha_3}(r))=\psi(yr).$

\begin{Claim}
For any $x\in X_3(F)$ and $y\in X_2(F)$ the spaces $\Hom_{V(\A)}(\Pi_3,\psi_x)$  and 
$\Hom_{N(\A)}(\Pi_3,\psi_y)$ are  one-dimensional. 
\end{Claim}

From the minimality of $\Pi_3$ it follows that the Fourier expansion of $\varphi \in \bar \theta(\Pi_3)$ reads
\begin{equation}\label{D3:D2:expansion}
\varphi=\sum_{x\in X_3(F)}\varphi^{V,\psi_x}+\varphi^V=
\sum_{x\in X_3(F)}\varphi^{V,\psi_x}+\sum_{y\in X_2(F)}\varphi^{N,\psi_y}+\varphi^N.
\end{equation}

For $\varphi=\bar\theta_{G,L}(f)$ one has 
$$\varphi^N(e)=\bar\mM_{w(21)}(f)+\bar\mM_{w(23)}(f)+ \bar E(f),$$ where 
$$\bar \mM_{w(21)}(f)=Res_{s=1}\mM_{w(21)}(s)(f)(e),\quad  \bar\mM_{w(23)}(f)=Res_{s=1}\mM_{w(23)}(s)(f)(e),$$
$$\bar E(f)=Res_{s=1}(\mM_{w(213)}(s)(f)(e)+\mM_{w(2132)}(s)(f)(e)).$$

\begin{Claim}
Let $\chi_1=w(21)^{-1}(\chi_{P,1})$ and $\chi_3=w(23)^{-1}(\chi_{P,1}).$
One has for $i=1,3$ 
$$\bar\mM_{2i}\in \Hom_{B(\A)}(\Ind^G_P(1), \chi_i),\quad  \bar\mM_{2i}(f^0)=\frac{R}{\xi(3)}.$$
\end{Claim}
The operators $\mM_{w(213)}(s)$ and $\mM_{w(2132)}(s)$ have both poles
of order $2$, but the sum has a simple pole.  
\begin{Lem}\label{D3:barE:kem}
The functional $\bar E\in \Hom_{N(\A)T(F)}(\Ind^G_P(1),\C)$ is uniquely determined by the condition:
For any $t=\Pi \alpha_i^\vee(r_i)$
$$\bar E(tf^0)=\frac{R^2}{\xi(2)\xi(3)}|r_1r_3|^{-1}|r_2|(2\ln|r_2|+A)$$
where $A=\left(\frac{\xi(s)}{\xi(s+1)}\right)'|_{s=0}$. 
\end{Lem}

\begin{proof}
We denote $\chi_0=w(213)^{-1}(\chi_{P,1})=w(2132)^{-1}(\chi_{P,1})=\chi_{-1,0,-1}.$

Let us rewrite $\mM_{w(2132)}(s)+\mM_{w(213)}(s)=(\mM_{w(2)}(s-1)+Id)\circ \mM_{213}(s)$.
The leading form of this operator is the composition of the leading terms of $\mM_{w(2)}(s-1)+Id$ and 
$\mM_{213}(s)$ at $s=1$.

The operator  $\mM_{213}(s): \Ind^G_P(s)\rightarrow \Ind^G_B \chi_{-s,s-1,-s}$ 
has a pole of order $2$ at $s=1$. The leading term defines a map 
$\mM^{-2}_{w(213)}:\Ind^G_P(1)\rightarrow \Ind^G_B \chi_0$, such that $\mM^{-2}_{w(213)}(f^0)=\frac{R^2}{\xi(2)\xi(3)}f_1^0$,
where $f^0, f^0_1$ are the normalized spherical sections in suitable representations.
Let $P_2=M_2 V_2$ be the standard parabolic subgroup, whose  Levi subgroup of $M_2$ is  generated by the roots $\pm \alpha_2$.
Rewrite 
$$\Ind^G_B\chi_{-s,s-1,-s}=
\Ind^G_{P_2} \left(|\omega_1|^{-s}\otimes |\omega_3|^{-s} \otimes \Ind^{M_2}_B |\omega_2|^{s-1}\right)$$
 The operator $\mM_{w(2)}(s-1)+Id$ acts on the inner induction and its leading term  
 is computed in Lemma \ref{M:der}. Hence $\bar E(tf^0)$ equals 
$$ Res_{s=1}(\mM_{w(2)(s-1)}+Id)\circ \mM_{213}(s)(tf^0)=\frac{R^2}{\xi(2)\xi(3)}|r_1r_3|^{-1}\cdot|r_2|(2\ln|r_2|+A),$$
as required.    
\end{proof}

Obviously there exists a functional $E\in \Hom_{N(\A)T(F)}(\mS(X_3(\A)),\C)$
such that  $\bar E=E\circ r_{3},$ so that 
$$\bar E(t\cdot \phi^0)=R\cdot |r_1r_3|^{-1}|r_2|(2\ln|r_2|+A).$$

Our next goal is to define the space $\mS(X_2(\A)),$ with a  surjective boundary map 
$\mB\in \Hom_{P(\A)}(\mS(X_3(\A),\mS(X_2(\A))$.

\subsection{Local Theory}
\subsubsection{Jacquet modules $(\Pi_{3,v})_N$ and $(\Pi_{3,v})_V$}
The structure of the Jacquet module $(\Pi_{3,v})_V$  is different the Jacquet modules of $\Pi_n$
for $n\ge 4.$ We shall describe it below. 

Denote
$$\chi_{1,v}=\chi_{(-1,-1,1)}= w(21)^{-1}(\chi_{P,1}),\quad 
\chi_{3,v}=\chi_{(1,-1,-1)}= w(23)^{-1}(\chi_{P,1}).$$
 $$\chi_{2,v}=\chi_{1,-2,1}=w(2)^{-1}(\chi_{P,1}), \quad \chi_{0,v}=\chi_{-1,0,-1}=
 w(213)^{-1}(\chi_{P,1})= w(2132)^{-1}(\chi_{P,1}).$$
Let $W$ be a two dimensional representation of  $T(F_v)$
with unique one-dimensional submodule and unique one-dimensional  quotient. 
The action of $T(F_v)$ on the submodule and the quotient are trivial. 
In some basis the action of $t=\Pi_{i=1}^3 \alpha_i^\vee(r_i)$ is given by 
the matrix $\left(\begin{array}{ll}
    1 & \ln|r_2|\\ 0 & 1  \end{array}\right)$

\begin{Claim}
\begin{enumerate}
\item  The representation
$\Ind^L_{B}(\chi_2)\simeq \Ind^{GL_2}_{B_2}\rho_{B_2}\otimes \Ind^{GL_2}_{B_2}(\rho_{B_2})|_L$ contains
unique irreducible subrepresentation $|\det|^{-2} (St\otimes St)|_L$ of $L(F_v).$
The quotient, denoted by $\Pi_{2,v}$ is a reducible representation of  $L(F_v)$ of length $3$.
\item
The Jacquet functor  $(\Pi_{3,v})_V$ fits in the non-split sequence of $L(F_v)$ modules:
$$0\rightarrow \C \rightarrow (\Pi_{3,v})_{V}\rightarrow \Pi_{2,v}\rightarrow 0,$$

\item
The Jacquet module $(\Pi_{3,v})_{N}$ is isomorphic to $\chi_{1,v}\oplus\chi_{3,v}\oplus \left(\chi_{0,v}\otimes W\right)$,

\end{enumerate}
\end{Claim}

\subsubsection{The space $\mS(X_2(F_v))$}
Let us describe a model of $\Pi_{2,v}$  in a space of functions $\mS(X_2(F_v))$ on $X_2(F_v).$

First we recall the Kirillov model 
of $\Ind^{GL_2}_{B_2}(s)$ on a space $\mS(F^\times)$, where  
Let $B_2=T_2N_2$ be the Borel subgroup of $GL_2$. 

For any $b\in F_v$ let $\psi_b$ be the  character of $F_v$ defined by 
$\psi_b(x)=\psi(bx).$ With fixed pinning it is considered as a character of $N_2(F_v).$

Define $J_b(s)\in \Hom_{N_2(F_v)}(\Ind^{GL_2}_{B_2}(|\rho_{B_2}|^s), \psi_b)$ by
$$J_b(f)=\int\limits_{N_2(F_v)} f(wn,s)\overline{\psi_b}(n)dn.$$
The integral converges absolutely for $Re(s)>0$ and is defined by meromorphic continuation for all 
$s$.  For $s=1$ the functional $J_0$ is $GL_2(F_v)$ invariant, whose kernel
is isomorphic to $St$, the Steinberg representation.

Consider a map 
$r:\Ind^L_B(\chi_2)\rightarrow \mS^\infty(X_2(F_v))$ defined by
$$r(f_1\otimes f_2)(b,0)=J_{b}(f_1) J_0(f_2),\quad
r(f_1\otimes f_2)(0,b)=J_0(f_1)J_{b}(f_2), $$
where $f_1\otimes f_2|_L \in \Ind^L_B(\chi_2).$
The kernel of this map is $|\det|^{-2} St\otimes St$ and  the image is denoted  by $\mS(X_2(F_v))$.
Hence $r$ serves an $L(F_v)$-isomorphism  $\Pi_{2,v}\simeq \mS(X_2(F_v)).$

\begin{Claim} 
\begin{enumerate}
\item
The normalized spherical function $\phi^0_{2,v}$ given by  
$$\phi^0_{2,v}(x)=\zeta_v(1)+\zeta_v(-1)\|ax\|_v,\quad  \|x\|\le |a|_v^{-1},$$
and zero otherwise,  generates the space $\mS(X_2(F_v)).$
The normalization is chosen such that $\phi_{2,v}^0(a^{-1},0)=\phi^0_{2,v}(0,a^{-1})=1$
\item
There exist functionals $\quad p_{i,v}\in \Hom(\mS(X_2(F_v)),\chi_{i,v})$ for $i=0,1,3$
such that  for any $\phi_v\in \mS(X_2(F_v))$ holds 
$$\phi_v(b,0)=p_{0,v}(\phi_v)+p_{1,v}(\phi_v)|b|_v, \quad \phi_v(0,b)=
p_{0,v}(\phi_v)+p_{3,v}(\phi_v)|b|_v$$
for $|b|_v\ll 1.$
In particular $p_{1,v}(\phi_{2,v}^0)=p_{3,v}(\phi_{2,v}^0)=\zeta_v(-1)|a|_v.$
\item The map $\phi_v\rightarrow (p_{0,v}(\phi_v),p_{1,v}(\phi_v), p_{3,v}(\phi_v))$
gives rise to the isomorphism $\mS(X_2(F_v))_{N_L}\simeq \chi_{0,v}\oplus \chi_{1,v}\oplus \chi_{3,v}.$
\end{enumerate}
\end{Claim}

\subsection{The adelic space $\mS(X_2(\A))$}
The global space $\mS(X_2(\A))$ is defined as the restricted tensor
product of local spaces with respect to the choice of $\phi^0_{2,v}$.
\begin{itemize}
    \item The global 
spherical function is denoted by $\phi^0_2$.
\item
The functional $\theta_{X_2}$ on $\mS(X_2(\A))$ is defined as in the introduction. 
\item
The  boundary map  is $\mB=\Pi_v\mB_v\in \Hom_{P(\A)}(\mS(X_3(\A),\mS(X_2(\A))$.
\item
The functionals $p_i\in \Hom_{B(\A)}(\mS(X_2(\A)),\chi_i)$ for $i=1,3$ 
are defined by $p_i(\phi^0_2)=\zeta(-1)|a|=\zeta(-1)q^{2-2g}.$
\end{itemize}

\begin{Lem}\label{D3:lem:X2:pi}
\begin{enumerate}
\item For $i=1,3$ one has $\bar\mM_{w(2i)}=|\Delta_F|^{\frac{1}{2}}\cdot p_i\circ \mB\circ r_3$.   
\item For any $y\in X_2(F)$ one has $\bar \theta_{G,L}(f)^{N,\psi_y}(f)(e)=
  |\Delta_F|^{\frac{1}{2}} \mB \circ r_3(f)(y).$
\end{enumerate}
\end{Lem}

\begin{proof} 
Since the space $\Hom_{B(\A)}(\mS(X_3(\A)),\chi_i)$ is one dimensional it is enough to check the 
equality on the spherical vector. 
$$p_i\circ \mB\circ r_3(f^0)=\zeta(-1)|a|\cdot \frac{R}{\xi(2)\xi(3)}=
\xi(-1)q^{-(1-g)}\cdot q^{2-2g} \cdot \frac{R}{\xi(2)\xi(3)}=$$
$$ q^{1-g}\cdot \frac{R}{\xi(3)}=
|\Delta_F|^{-\frac{1}{2}}\mM_{w(2i)}.$$
Here we have used  the functional equation $\xi(-1)=\xi(2).$

Let $y=(b,0)\in X_2(F)$. We write $N=N'\cdot N_{\alpha_1}.$
Then $\bar\theta_{G,L}(f^0)^{N,\psi_y}(f^0)(e)$ is computed in stages as 
$(\bar\theta_{G,L}(f^0)^{N'})^{N_{\alpha_1},\psi_y}(e).$ Only one element $w(23)\in W_L\backslash W/W_{\alpha_1}$
has a non-zero contribution which equals 
$$Res_{s=1} \frac{\xi(s)}{\xi(s+2)}\mE^{GL_2}_B(f^0,s)^{N,\psi_y}=q^{-(1-g)}\frac{R}{\xi(3)\xi(2)}\sigma(ay,-1)$$
On the other hand 
$$\mB\circ r_3(f^0)(y)=\frac{R}{\xi(3)\xi(2)}\phi^0_2(y)=\frac{R}{\xi(3)\xi(2)}\phi^0_2(y)\sigma(ay,-1).$$
So 
$$\bar \theta_{G,L}(f)^{N,\psi_y}(f^0)(e)=  q ^{-(1-g)} \mB \circ r_3(f^0)(y)=
|\Delta_F|^{\frac{1}{2}}\mB \circ r_3(f^0)(y).$$

    \end{proof}

Now we are ready to formulate the analogue of Theorem \ref{main:thm} for the pair $(D_3,D_2).$ 
\begin{Thm} The following decomposition holds for the functional $\theta_{D_3,D_2}$ on $\mS(X_3(\A))$
  $$\theta_{D_3,D_2}=\theta_{X_3}+|\Delta_F|^{\frac{1}{2}}\theta_{X_2}\circ \mB+
  |\Delta_F|^{\frac{1}{2}}(p_1+p_3)\circ \mB+ E$$
    \end{Thm}
\begin{proof} The Theorem follows from the Fourier expansion \ref{D3:D2:expansion}
and Lemmas \ref{D3:barE:kem},\ref{D3:lem:X2:pi}.
\end{proof}

\section{Transition between models}\label{sec:transition}
In Theorem \ref{main:thm} there is a boundary term $\theta_{G_1,M_1}\circ \mB$.
If the pair $(G_1,M_1)$ is not admissible, Theorem \ref{main:thm} can not
be applied to it. Hence it is desirable to replace the functional
$\theta_{G_1,M_1}$ by another automorphic functional on $\Pi_1$
associated to an admissible pair.
In this section we explore relations between functionals
$\theta_{G,L_1}$ and $\theta_{G,L_2}$ for two weakly admissible pairs $(G,L_1)$ and $(G,L_2)$
with parameters $(s_1,d_1)$ and $(s_2,d_2).$
They give rise to the models $\mS(X_1(\A))$ and $\mS(X_2(\A))$
of the same 
minimal representation $\Pi$ of $G(\A)$.
Let $\phi^0_1,\phi^0_2$ be the normalized
spherical functions in each model. There is a $G(\A)$-equivariant isomorphism 
\beq \label{transition:def}
T:\mS(X_1(\A))\rightarrow \mS(X_2(\A)),\quad T(\phi_1^0)=\phi^0_2 
\eeq
The functionals $\theta_{G,L_1}$ and $\theta_{G,L_2}\circ T$ on $\mS(X_1(\A))$ are proportional. 
Let us compute the constant of proportionality. 
\begin{Thm} \label{theta:transition}
One has $\theta_{G,L_1}=\kappa\cdot \theta_{G,L_2}\circ T$, where 
\begin{itemize}
    \item $G=D_n, L_1=A_{n-1}, L_2=D_{n-1}, \quad \kappa=|\Delta_F|^{\frac{n-4}{2}}.$ 
    \item $G=E_6, L_1=A_5, L_2=D_5 ,\quad \kappa=|\Delta_F|^{\frac{1}{2}}.$
    \item $G=E_7, L_1=D_6, L_2=E_6 , \quad\kappa =|\Delta_F|^{\frac{1}{2}}.$
\end{itemize}

\end{Thm}

\begin{proof}
The result follows from an unexpected identity between values 
of different degenerate Whittaker  functionals applied
to a spherical automorphic
function $\varphi^0$ in the minimal representation.  

For the triple $(G,L_1,L_2)$ let $\beta_i$ be the simple root defining $L_i.$
Write $N=N_i\cdot N_{\beta_i}$. Let $\psi_i$ be the character of $N(F)\backslash N(\A)$
such that $\psi_i|_{N_{\beta_i}}=\psi$ and $\psi|_{N_i}=1.$
The degenerate spherical Whittaker coefficient is defined by 
$$W_i(\varphi^0)=\integral{N} \varphi^0(n)\overline{\psi_i(n)} \,dn,$$ 
where $\varphi^0$ is a spherical vector in $\bar\theta(\Pi).$

\begin{Prop}\label{whittaker:transition} For the pairs $(G,L_1)$ and
  $(G,L_2)$ as above one has 
$$\frac{W_1(\varphi^0)}{\sigma(a,d_1)}=\kappa\cdot 
\frac{W_{2}(\varphi^0)}{\sigma(a,d_{2})},$$
where $\kappa$ is defined as in Theorem \ref{theta:transition}.
\end{Prop}

Let us derive Theorem \ref{theta:transition}, from Proposition \ref{whittaker:transition} 

The spherical functions $f^0_i\in \Ind^{G}_{P'_i}(s_i)$ are normalized by 
$f^0_i(e)=1$ for $i=1,2.$
The functions $\bar\theta_{G,L}(f^0_1)(g)$ and $\bar \theta_{G,L_2}(f_2^0)(g)$
are spherical functions in the irreducible automorphic representation
$\bar\theta(\Pi)$ and hence are proportional. 
Write  
$$K\cdot \frac{\bar\theta_{G,L}(f^0_1)(g)}{c_{G,L_1}}=
\frac{\bar \theta_{G,L_2}(f_2^0)(g)}{c_{G,L_2}}.$$
One has 
$$K=K\cdot \frac{\bar\theta_{G,L_1}(f_1^0)^{N,\psi_1}(e)}{c_{G,L_1}\sigma(a,d_1)}=
\frac{\bar\theta_{G,L_2}(f_2^0)^{N,\psi_1}(e)}{c_{G,L_2}\sigma(a,d_1)}=\kappa\cdot \frac{\bar\theta_{G,L_2}(f_2^0)^{N,\psi_2}(e)}{c_{G,L_2}\sigma(a,d_2)}=\kappa.$$
The third equality follows from Proposition \ref{whittaker:transition}, applied
to $\varphi^0=\theta_{G,L_2}(f^0_2)$. Hence $K=\kappa$. By definition of 
the functionals $\theta_{G,L_i}$ one has 
$$\theta_{G,L_1}(\phi_1^0)=
\frac{\bar\theta_{G,L_1}(f^0_1)(e)}{c_{G,L_1}}=\kappa\cdot 
\frac{\bar\theta_{G,L_2}(f^0_2)(e)}{c_{G,L_2}}=\kappa\cdot \theta_{G,L_2}(\phi^0_2)$$
and so $\theta_{G,L_1}=\kappa \cdot \theta_{G,L_2}\circ T$ as required.

\end{proof}

It remains to prove Proposition \ref{whittaker:transition}.
\begin{proof} The argument is of inductive nature. 
  The proof for $D_n$  is by induction on $n$.
  The proof for $E_6$ relies on the result for $D_5$
  and the proof for $E_7$ relies on the result for $E_6$. 
In all the cases we use the identity (\ref{second:identity}).
  
\begin{enumerate}
    \item
Let $G_n$ be a simply-connected group of type $D_n$ and the maximal parabolic
$P_1=L_1V_1$ such that $[L_1,L_1]$ has type $D_{n-1}$.
Let $N_n$ be the unipotent radical of the Borel subgroup in $G_n$. 
We have $d_1=n-3$ and $d_2=1$.
Let $\varphi^0=\bar\theta_{D_n,D_{n-1}}(f^0).$
By Theorem \ref{spherical:Fourier}
$$W_1(\varphi^0)=c_{D_n,D_{n-1}}\cdot \sigma(a,n-3).$$
The goal is to show that
$$W_2(\varphi^0)=q^{(1-g)(n-4)}\cdot c_{D_n,D_{n-1}}\cdot \sigma(a,1).$$

The proof is by induction on $n$. For $n=4$ one has $d_1=d_2=1.$
The representation $\Pi$ can be extended to the representation
$\Pi$ of $L(\A)\rtimes S_3$.
and the group $S_3$ preserves the spherical function $\varphi^0$
and the unipotent radical $N$.
The characters $\psi_1$ and $\psi_2$ belong to the same $S_3$
orbit and hence 
$W_1(\varphi^0)=W_2(\varphi^0)$, as required.

One has
$$\bar\theta_{D_n,D_{n-1}}(f^0)^{N_n,\psi_2}=
(\bar\theta_{D_n,D_{n-1}}(f^0)^{V_1})^{N_{n-1},\psi_{2}}.$$
Let us first compute the inner integral
$$\bar\theta_{D_n,D_{n-1}}(f^0)^{V_1}(g)=Res_{s=s_0} c_{v_0}(s)f^0(g)+
c_{v_1}(1)\cdot \bar\theta_{D_{n-1},D_{n-2}}(f^0)$$

The first summand is a constant function of $N_{n-1}$,
so its  Whittaker coefficient is zero. By induction assumption
$$\bar\theta_{D_{n-1},D_{n-2}}(f^0)^{N_{n-1},\psi_2}(e)=
q^{(1-g)(n-5)} c_{D_{n-1},D_{n-2}}\sigma(a,1)$$
and by equation \ref{second:identity}
$c_{v_1}(1)\cdot c_{D_{n-1},D_{n-2}} =  q^{(1-g)}  \cdot c_{D_n,D_{n-1}}$.
So $$W_2(\varphi^0)=q^{(1-g)(n-4)}\cdot c_{D_n,D_{n-1}}\sigma(a,1)$$
and we are done.

\item
Let $G$ be of type $E_6$ and  $P_1=L_1V_1,$ where $L_1$ is of type $D_5$.
We pick the spherical function $\varphi^0=\bar\theta_{E_6,D_5}(f^0)$.
By Theorem \ref{spherical:Fourier} one has 
$W_1(\varphi^0)=c_{E_6,D_5}\cdot \sigma(a,2)$.

One has $N=V_1\cdot N',$ where  $N'$ is  the unipotent radical of Borel subgroup of $L_1$.  Let us compute $(N,\psi_2)$ coefficient in stages. 
$$W_2(\varphi^0)=(\bar\theta_{E_6,D_5}(f^0)^{V_1})^{N',\psi_2}(e)=
c_{v_1}(3)\bar\theta_{D_5,A_4}(f^0)^{N',\psi_2}(e)=$$
$$=c_{v_1}(3)c_{D_5,A_4}\cdot \sigma(a, 1)=
q^{g-1}\cdot c_{E_6,D_5}\cdot \sigma(a,1).$$
\item
Finally, let $L=E_7$ and $P_1=L_1V_1$ when $L_1$ is of type $E_6$.
For $\varphi^0=\bar\theta_{E_7,E_6}(f_1^0)$
one has $W_1(\varphi^0)=c_{E_7,E_6}\cdot \sigma(a,3)$.
Let us decompose $N=V_1\cdot N',$ where $N'$ is  the unipotent radical of Borel subgroup 
of $L_1$.

$$W_2(\varphi^0)=(\bar\theta_{E_7,E_6}(f^0)^{V_1})^{N',\psi_2}(e)=
c_{v_1}(5)\bar\theta_{E_6,D_5}(f^0)^{N',\psi_2}(e)=$$
$$  c_{v_1}(5)c_{E_6,D_5} \sigma(a,2)=q^{g-1} c_{E_7,E_6}\sigma(a,2),$$
as required. Here we applied again the equation \ref{second:identity}
\end{enumerate}

\end{proof}

\section{Automorphic functionals for $E_6,E_7,E_8$}\label{sec:autom:func}
We can apply Theorems \ref{main:thm} and
\ref{theta:transition} for writing the expansion 
of the automorphic functionals on the groups $E_6,E_7,E_8$.
\subsection{$G=E_6$}
Let $X,X_1,X_2$ be the models for the pairs 
$(E_6,D_5)$, $(D_5,A_4)$ and $(D_5,D_4)$ respectively.

Let  $\mB_i:\mS(X(\A))\rightarrow \mS(X_i(\A))$ be the boundary maps for
$i=1,2$ and $T:\mS(X_1(\A))\rightarrow \mS(X_2(\A))$ be the transition map so that 
$\mB_2=\mB_1\circ T.$

Applying Theorem \ref{main:thm} we have 
$$\theta_{E_6,D_5}=\theta_X+p_X+|\Delta_F|^{\frac{1}{2}} \theta_{D_5,A_4}\circ \mB_1.$$

The pair $(D_5,A_4)$ is not admissible.
Applying further Theorem \ref{theta:transition} for the pairs $(D_5,A_4)$ and $(D_5,D_4)$
we arrive to 
\beq \label{E6:automorpic:decomp}
\theta_{E_6,D_5}=\theta_X+p_X+ |\Delta_F| \theta_{D_5,D_4}\circ \mB_2.
\eeq
Since the pair $(D_5,D_4)$ is admissible we can apply Theorem \ref{Dn:automorphic:decomp}
for the full expansion.

\subsection{$G=E_7$}
Applying  Theorem \ref{main:thm} to the functional
$\theta_{E_7,E_6}$ we get
$$\theta_{E_7,E_6}=\theta_X+p_X+|\Delta_F|^{1/2}\theta_{E_6,D_5}\circ \mB$$
and the expansion of $\theta_{E_6,D_5}$ has been discussed above.

Alternatively, we can address the functional $\theta_{E_7,D_6}$.
Let $X,X_1,X_2$ be the models for the pair 
$(E_7,D_6)$, $(D_6,A_5)$ and $(D_6,D_5)$ respectively.
Let  $\mB_i:\mS(X(\A))\rightarrow \mS(X_i(\A))$ be the boundary maps for $i=1,2$
and $T:\mS(X_1(\A))\rightarrow \mS(X_2(\A))$ be the transition map so that 
$\mB_2=\mB_1\circ T.$

Applying Theorems \ref{main:thm} and \ref{theta:transition} we see 
 the boundary term of $\theta_{E_7,D_6}$ can be rewritten
$$p_Y\circ \iota+|\Delta_F|^{\frac{1}{2}}\theta_{D_6,A_5}\circ \mB_1=p_Y\circ \iota+
|\Delta_F|^{\frac{3}{2}}\cdot \theta_{D_6,D_5}\circ \mB_2.$$
Since the pair $(D_6,D_5)$ is admissible, the expansion
can be continued by Theorem \ref{Dn:automorphic:decomp}.

\subsection{$G=E_8$}
The group $E_8$ does not have a Siegel parabolic subgroup, and
the only functional to consider is $\theta_{E_8,E_7}$. 
The boundary term of the expansion is expressed by Theorem \ref{main:thm}
in terms of $\theta_{E_7,E_6},$ that has been discussed above. 

\bibliographystyle{alpha}
 \bibliography{bib}

\end{document}